%% file: minimal.tex
\documentclass[12pt]{amsart}
\usepackage[english]{babel}
\usepackage[latin1]{inputenc}

\usepackage{amsmath}
\usepackage{amsfonts}
\usepackage{amssymb}
\usepackage{amsthm}
\usepackage{graphicx,epsfig}
\usepackage{amscd}

\DeclareMathOperator{\re}{Re}
\DeclareMathOperator{\im}{Im}

\DeclareMathOperator{\hg}{F}
\DeclareMathOperator{\tr}{tr}

\newcommand{\R}{\mathbb{R}}
\newcommand{\h}{\mathbb{H}}
\newcommand{\Z}{\mathbb{Z}}
\newcommand{\N}{\mathbb{N}}
\newcommand{\C}{\mathbb{C}}
\newcommand{\s}{\mathbb{S}}

\newcommand{\rmS}{\mathrm{S}}
\newcommand{\rmd}{\mathrm{d}}
\newcommand{\rmO}{\mathrm{O}}

\newcommand{\rmD}{\mathrm{D}}

\newcommand{\rmI}{\mathrm{I}}
\newcommand{\rmII}{\mathrm{II}}

\newcommand{\cS}{{\mathcal S}}
\newcommand{\cK}{{\mathcal K}}

\newcommand{\cT}{{\mathcal T}}

\newcommand{\cD}{{\mathcal D}}

\newcommand{\cO}{{\mathcal O}}
\newcommand{\cL}{{\mathcal L}}
\newcommand{\cP}{{\mathcal P}}
\newcommand{\cI}{{\mathcal I}}

\begin{document}

\newtheorem{thm}{Theorem}
\newtheorem{defn}[thm]{Definition}
\newtheorem{cor}[thm]{Corollary}
\newtheorem{prop}[thm]{Proposition}
\newtheorem{app}[thm]{Application}
\newtheorem{lemma}[thm]{Lemma}
\newtheorem{rem}[thm]{Remark}
\newtheorem{ex}[thm]{Example}
\newtheorem{notation}[thm]{Notations}
\newtheorem{hypothesis}[thm]{Hypothesis}

\title[Minimal disks bounded by three straight lines]
{Minimal disks bounded by three straight lines in Euclidean
space and trinoids in hyperbolic space} 
\author{Beno\^\i t Daniel}
% \date{}

\begin{abstract}
Following Riemann's idea, we prove the existence of a minimal disk
in Euclidean space bounded by three lines in generic position and with
three helicoidal ends of angles less than $\pi$. In the case of
general angles, we prove that there exist at most four such minimal
disks, we give a sufficient condition of existence in terms of 
a system of three equations of degree $2$, and we give explicit
formulas for the Weierstrass data in terms of hypergeometric
functions. Finally, we construct constant-mean-curvature-one trinoids
in hyperbolic space by the method of the conjugate cousin immersion. 
\end{abstract}

\subjclass{Primary: 53A10. Secondary: 53C42, 53A35, 30F45}
\keywords{Minimal surfaces, Euclidean space, Weierstrass representation, Bryant
surfaces, hyperbolic space, constant mean curvature}

\maketitle

\section{Introduction}

In this paper we investigate minimal disks in Euclidean space $\R^3$ bounded by
three straight lines in generic position ({\it i.e.} the lines do not intersect
one another and do not lie in parallel planes, in particular they are
not pairwise parallel).
This problem was investigated by Riemann in his posthumous
memoir \cite{riemann}. He actually introduced the spinor representation,
the Gauss map and the Hopf differential of minimal surfaces in
Euclidean space. 
He investigated the case of minimal
surfaces bounded by a contour composed of pieces of straight lines (possibly
going to infinity). He studied more precisely the cases where the
contour is composed of $2$, $3$ or $4$ lines.

However, his study was not complete and sometimes not precise; in
particular, he did not deal with questions of orientations. The first
aim of this paper is to complete Riemann's study of minimal surfaces
bounded by three straight lines. More precisely we will investigate
minimal immersions $x$ from $\Sigma=\{z\in\C|\im z\geqslant
0\}\setminus\{0,1\}$ into $\R^3$ mapping $(-\infty,0)$, $(0,1)$ and
$(1,\infty)$ onto three lines (in generic position) and having
helicoidal ends (in the sense explained in section
\ref{helicoidalends}) at $0$, $1$ and $\infty$. The method is to study
the Weierstrass data of $x$ and to use the spinor representation.

We first prove that there exist at most four minimal immersions
bounded by three given lines (in generic position) with helicoidal
ends of given parameters (theorem \ref{atmostfour}). To do this, we
use a result by Riemann: he proved that the spinor data
satisfy a differential equation involving the Schwarzian derivative of
the Gauss map. This is a second order equation with five regular singularities.
Studying the behaviour of the Schwarzian derivative at these singular
points, we prove that the Schwarzian derivative only depends on two
parameters that are related by two polynomial equations of degree $2$,
and thus that there are at most four possiblities for the Schwarzian
derivative of the Gauss map.

We next study the explicit immersion given by Riemann in his memoir
\cite{riemann}: he introduced spinors in terms of hypergeometric
functions, but he did not check that they actually gave a minimal
immersion bounded by given lines. We establish this fact in proposition
\ref{threelines}. More precisely, given three lines in generic
position, denoting by $A$, $B$, $C$ the distances between the lines, and
$\pi\alpha$, $\pi\beta$, $\pi\gamma$ the angles of the ends (with
signs as explained in section \ref{helicoidalends}),
we prove that to each real solution
$(p,q,r)$ (up to a sign) of the system
\begin{eqnarray*}
\left\{\begin{array}{ccc}
p^2-\alpha^2(p+q+r)^2 & = & \varepsilon\frac{A\alpha}{2\pi} \\
q^2-\beta^2(p+q+r)^2 & = & \varepsilon\frac{B\beta}{2\pi} \\
r^2-\gamma^2(p+q+r)^2 & = & \varepsilon\frac{C\gamma}{2\pi}
\end{array}\right.
\end{eqnarray*}
where $\varepsilon\in\{1,-1\}$ (depending on the geometric configuration of
the lines and on the angles) corresponds a
minimal immersion (with possibly a singular point when the Hopf
differential has a double zero, which is a non-generic situation)
bounded by these lines and with helicoidal ends of
parameters $(A,\alpha)$, $(B,\beta)$, $(C,\gamma)$ (theorem
\ref{theoremsystem}). In particular we prove the existence of at least
one minimal immersion in the case where the angles are less than $\pi$
(corollary \ref{existencesmallangles}). 

However we do not know if we obtain all the solutions in this way.

\begin{figure}[htbp]
\begin{center}
\includegraphics{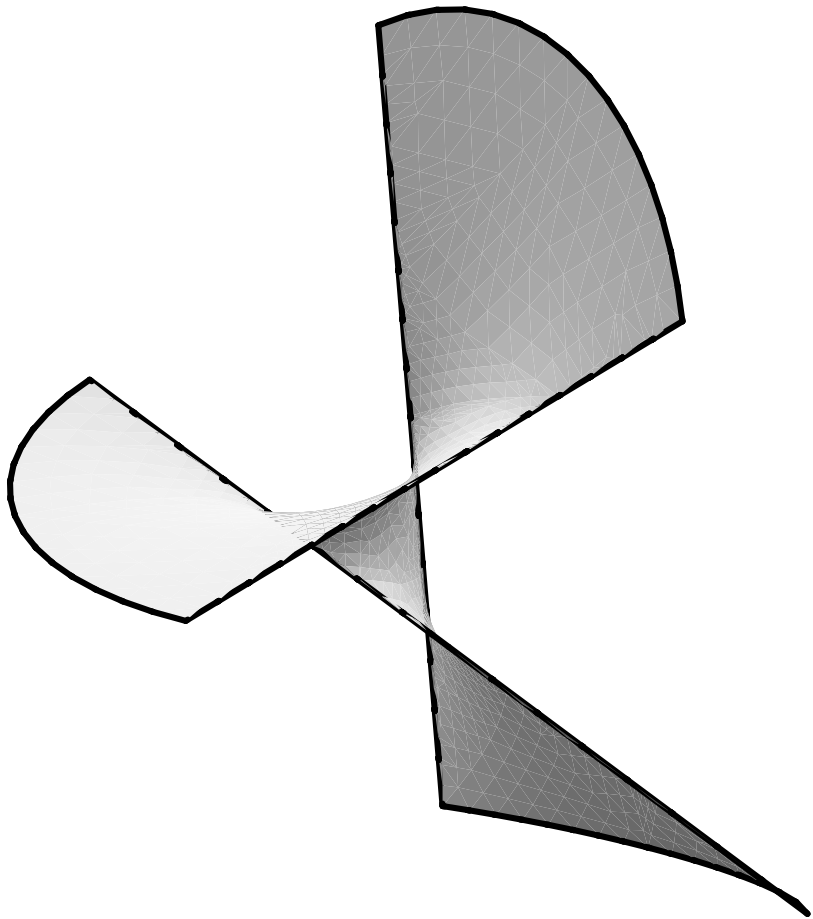}
\caption{A minimal surface with three helicoidal ends.}
\label{figure3lines}
\end{center}
\end{figure}

Figure \ref{figure3lines} is a
picture of a minimal surface bounded by three lines with helicoidal
ends, drawn with the software ``Evolver''.

In the last part of this paper, we construct
constant-mean-curvature-one (CMC-$1$) trinoids in hyperbolic space $\h^3$
applying the conjugate cousin method to minimal disks bounded by three
straight lines in $\R^3$. 
CMC-$1$ surfaces in $\h^3$ are called Bryant surfaces.
Bryant proved in \cite{bryant} that they are closely related to
minimal surfaces in $\R^3$ (see also \cite{umehara} and
\cite{rosenberg}); in particular there exists a representation in
terms of holomorphic data analogous to Weierstrass representation.

Irreducible trinoids in $\h^3$ were first classified by
Umehara and Yamada in \cite{metrics}, and then by Bobenko,
Pavlyukevich and Springborn in \cite{bobenko}, using different
techniques: their method has some similarities with that
used in this paper to find minimal surfaces bounded by three lines in
$\R^3$ (they use a spinor representation for Bryant surfaces and they
obtain explicit formulas in terms of hypergeometric 
functions). The technique of the conjugate cousin immersion was used
by Karcher in \cite{karcher} to construct trinoids with dihedral
symmetry. Here we 
use this technique to construct general irreducible trinoids with a
symmetry plane (actually every irreducible trinoid has a symmetry
plane by the classification of \cite{metrics}). We
also prove that the asymptotic boundary points
of the ends of these trinoids are distinct (except in
exceptionnal cases) (theorem \ref{theoremtrinoids}). Finally we give
examples of minimal disks bounded by three lines whose conjugate
cousins are invariant by some parabolic isometries.

{\bf Acknowledgements.} The author is grateful to Pascal Collin for
his explainations on Riemann's memoir \cite{riemann}, and to
Harold Rosenberg for submitting this problem and for discussions about
this paper.

\section{Preliminaries}

In all this paper, we will set $\bar\C=\C\cup\{\infty\}$,
$\bar\R=\R\cup\{\infty\}$,
$$\Sigma=\{z\in\C|\im z\geqslant 0\}\setminus\{0,1\},\quad
\Sigma_0=\{z\in\Sigma||z|<1\},$$
$$\Sigma_1=\{z\in\Sigma||z-1|<1\},\quad
\Sigma_{\infty}=\{z\in\Sigma||z|>1\}.$$
The canonical scalar product of $\R^3$ is denoted by $<\cdot,\cdot>$.
If $D$ is a straight line in $\R^3$, then $D^\perp$ denotes the set of
unit vectors that are orthogonal to $D$. The canonical basis of $\R^3$
will be denoted by $(\vec{e}_1,\vec{e}_2,\vec{e}_3)$:
$$\vec{e}_1=(1,0,0),\quad\vec{e}_2=(0,1,0),\quad\vec{e}_3=(0,0,1).$$
We define the logarithm and non-integer powers on $\Sigma$ in
the following 
way. Let $\kappa\in\R$. If $z=\rho e^{i\theta}\in\Sigma$ with $\rho>0$
and $\theta\in[0,\pi]$, then $\ln z=\ln\rho+i\theta$ and $z^\kappa=\rho^\kappa
e^{i\kappa\theta}$. If $z\in\Sigma$ and $z-1=\rho e^{i\theta}$ with $\rho>0$
and $\theta\in[0,\pi]$, then $(z-1)^\kappa=\rho^\kappa
e^{i\kappa\theta}$ and 
$(1-z)^\kappa=\rho^\kappa e^{i\kappa(\theta-\pi)}
=e^{-i\pi\kappa}(z-1)^\kappa$ (this convention is chosen in order that
$(1-z)^\kappa$ is real when $z$ is real and less than $1$).

Finally, $\cD$ denotes the set of the triples of straight lines in
$\R^3$ that are neither pairwise concurrent, neither pairwise
parallel, nor lying in parallel planes, modulo direct isometries of $\R^3$.

\subsection{Weierstrass representation} \label{weierstrass}

In this section we recall basic facts about Weierstrass representation
and we introduce some notations.

Let $\cS$ be a Riemann surface with boundary, and
$x=(x_1,x_2,x_3):\cS\to\R^3$
a conformal minimal immersion. Then we have
$$x(z)=x(z_0)+\re\int_{z_0}^z\left((1-g^2),i(1+g^2),2g\right)\omega$$
where $z_0$ is a fixed point in $\cS$ and $(g,\omega)$ the
Weierstrass data of $x$: $g$ is a meromorphic function
on $\cS$ and $\omega$ a holomorphic $1$-form on $\cS$. The poles of $g$
are the zeros of $\omega$, and $z$ is a pole of $g$ of order $k$ if
and only if $z$ is a zero of $\omega$ of order $2k$. Conversely, if
$g$ and $\omega$ satisfy this condition, then they define a minimal immersion.

We define $X=(X_1,X_2,X_3):\cS\to\C^3$ by
$$X(z)=x(z_0)+\int_{z_0}^z\left((1-g^2),i(1+g^2),2g\right)\omega.$$
We have $$\rmd x_1+i\rmd x_2=\bar{\omega}-g^2\omega,\quad
\rmd x_3=\re(2g\omega).$$
The Gauss map of $x$ is 
$$N=\left(\frac{2g}{|g|^2+1},\frac{|g|^2-1}{|g|^2+1}\right).$$
The orientations induced on $x(\cS)$ by the Gauss map $N$ and the
immersion $x$ are compatible. The function $g$ is the composition of
the Gauss map and the stereographic projection with respect to the
north pole of the sphere. If $D$ is a line parallel to the vector
$(\cos(\pi\alpha),\sin(\pi\alpha),0)$ for $\alpha\in\R$, then the circle
$D^\perp$ corresponds to $g\in ie^{i\pi\alpha}\bar\R$.

If $f$ is a meromorphic function on an open set $U\subset\cS$, we
define its Schwarzian derivative with respect to a local conformal
coordinate $z$ by
$$\rmS_zf=\left(\left(\frac{f''}{f'}\right)'
-\frac{1}{2}\left(\frac{f''}{f'}\right)^2\right)\rmd z^2.$$
If $\zeta$ is another local conformal coordinate, then
$\rmS_zf=\rmS_\zeta f+\rmS_z\zeta$. If $f$ is regular at a point
$z_1$, then $\rmS_zf$ is holomorphic at $z_1$; if $f$ has a branch
point of order $j-1$ at $z_1$ with $j\geqslant 2$, then $\rmS_zf$ has a
pole of order $2$ at $z_1$, and its coefficient of order $-2$ is equal
to $\frac{1-j^2}{2}$.

The Hopf differential is the holomorphic $2$-form on $\cS$ defined by
$$Q=\omega\rmd g=\frac{1}{2}\rmd X_3\frac{\rmd g}{g}.$$
The forms $Q$ and $\rmS_zg$ are invariant by a direct isometry of
$\R^3$. The first and second fondamental forms of the surface given by
$$\rmI=(1+|g|^2)^2|\omega|^2,\quad\rmII=-2\re Q.$$

\begin{prop} \label{zerohopf}
Let $z\in\cS$ and $k\in\N^*$. Then $z$ is a zero of $Q$ of order
$k$ if and only if $z$ is a branch point of $g$ of order $k$. This
happens if and only if one of the following conditions holds:
\begin{itemize}
\item $z$ is a zero of $g$ of order $k+1$,
\item $z$ is a pole of $g$ of order $k+1$,
\item $z$ is a zero of $\frac{\rmd g}{g}$ of order $k$.
\end{itemize}
\end{prop}

\begin{proof} If $z$ is a zero of $g$ of order $d\in\N^*$, then it is
not a zero of $\omega$ and it is a simple pole of $\frac{\rmd
g}{g}$. Consequently it is a zero of $Q$ of order $k$ if and only if $k=d-1$.

If $z$ is a pole of $g$ of order $d\in\N^*$, then it is a zero of
$\omega$ of order $2d$ and it is a simple pole of $\frac{\rmd
g}{g}$. Consequently it is a zero of $Q$ of order $k$ if and only if $k=d-1$.

If $z$ is neither a zero nor a pole of $g$, then it is not a zero of
$\omega$, and consequently it is a zero of $Q$ of order $k$ if and
only if it is a zero of $\frac{\rmd g}{g}$ of order $k$.
\end{proof}

We recall the spinor representation of a minimal surface:
$$x(z)=x(z_0)+\re\int_{z_0}^z
\left(\xi_1^2-\xi_2^2,i(\xi_1^2+\xi_2^2),2\xi_1\xi_2\right)$$
where $\xi_2$ and $\xi_2$ are holomorphic sections of a spin
structure (see \cite{kusner} for more information). These spinors satisfy
$$g=\frac{\xi_2}{\xi_1},\quad\omega=\xi_1^2.$$ Two such holomorphic
spinors define a minimal immersion if and only if they do not have
common zeros (if they have a common zero, then the map $x$ is not an
immersion at this point). We will call $(\xi_1,\xi_2)$ the spinor data
of $x$.

\subsection{Helicoidal ends} \label{helicoidalends}

Most of the results of this section are contained in Riemann's memoir
\cite{riemann}. Here we prove these results using modern formalism,
and we give precise definitions for helicoidal ends and the signs of
distances and angles. 

\begin{defn} \label{associatedvector}
Let $D_1$ and $D_2$ be two nonparallel and nonconcurrent oriented
straight lines, and let 
$\vec{u}_1$ and $\vec{u}_2$ be unit vectors inducing the orientations
of $D_1$ and $D_2$. Then the unit vector
$$\vec{v}=\frac{\vec{u}_1\times(-\vec{u}_2)}
{||\vec{u}_1\times\vec{u}_2||}$$ is called the vector associated to
the couple $(D_1,D_2)$ of oriented straight lines.
\end{defn}

\begin{defn} \label{signeddistance}
The signed distance of $D_1$ and $D_2$ is the number
$\rmD(D_1,D_2)=<\overrightarrow{p_1p_2},\vec{v}>$ 
where $p_1\in D_1$ and $p_2\in D_2$ (this number does
not depend on the choices of $p_1$ and $p_2$, and $|\rmD(D_1,D_2)|$ is
the distance between $D_1$ and $D_2$).
\end{defn}

\begin{defn} \label{helicoidalend} 
Let $U$ be a neighbourhood of $0$ in $\C$ that is symmetric with
respect to the real axis ({\it i.e.} $z\in U\iff\bar z\in U$),
$\Omega=\Sigma\cap U$,
$\Omega_1=\Omega\cap(-\infty,0)$, $\Omega_2=\Omega\cap(0,+\infty)$ and
$x:\Omega\to\R^3$ be a conformal minimal immersion that is complete at
$0$. Let $D_1$ and $D_2$ be two nonparallel and nonconcurrent oriented
straight lines, let $\vec{u}_1$ and $\vec{u}_2$ be unit vectors
inducing the orientations of $D_1$ and $D_2$, and let $\vec{v}$ be the
vector associated to $(D_1,D_2)$.

We say that the immersion $x$ has an end (at $z=0$) bounded
by the couple of oriented lines $(D_1,D_2)$ if
\begin{enumerate}
\item[1.] the immersion $x$ maps $\Omega_1$ to a part of $D_1$ and
$<x(z),\vec{u}_1>\to+\infty$ when $z\to 0$ with $z$ real and negative,
\item[2.] the immersion $x$ maps $\Omega_2$ to a part of $D_2$ and
$<x(z),\vec{u}_2>\to-\infty$ when $z\to 0$ with $z$ real and positive.
\end{enumerate}

We say that the immersion $x$ has a helicoidal end (at $z=0$) bounded
by the couple of oriented lines $(D_1,D_2)$ if moreover the two following
conditions are satisfied:
\begin{enumerate}
\item[3.] the Gauss map of $x$ has a limit when $z\to 0$,
\item[4.] the quantity $<x(z),\vec{v}>$ is bounded when $z\to 0$.
\end{enumerate}
\end{defn}

It will follow from the proof of lemma \ref{parameters} that a
helicoidal end is actually asymptotic to a helicoid.

\begin{lemma} \label{gaussend}
Assume that $x$ has a helicoidal end bounded by $(D_1,D_2)$. Let
$N$ be the Gauss map of $x$. Then the limit point of $N$ at $0$ is
$N(0)=\vec{v}$ or $N(0)=-\vec{v}$.
\end{lemma}

\begin{proof}
We have $N(z)\in D_1^\perp$ if $z\in\Omega_1$ and $N(z)\in D_2^\perp$
if $z\in\Omega_2$, so $N(0)\in D_1^\perp\cap D_2^\perp$.
\end{proof}

\begin{lemma} \label{hopfschwarziansymmetry}
Assume that $x$ has an end bounded by $(D_1,D_2)$. Let
$(g,\omega)$ be the Weierstrass data of $x$, and $Q$ its Hopf
differential. Then $Q$ extends to a holomorphic $2$-form on
$U\setminus\{0\}$ and $\rmS_z g$ extends to a meromorphic $2$-form on
$U\setminus\{0\}$.
\end{lemma}

\begin{proof}
Since a direct isometry of $\R^3$ does not change $Q$ and $\rmS_z g$, we
can assume that $D_2$ is the $x_1$-axis. Then $x_3$ is
constant on $\Omega_2$, so $x_3'=0$ on $\Omega_2$, and thus
$X_1'(z)\in i\R$ for $z\in\Omega_2$. And the Gauss map $N$ is normal
to the $x_1$-axis on $\Omega_2$, so $g(\Omega_2)\subset i\bar\R$, and
thus $\frac{g'}{g}(z)\in\bar\R$ for $z\in\Omega_2$. Thus
$\frac{Q}{\rmd z^2}=\frac{1}{2}X_3'\frac{g'}{g}$ is purely imaginary
on $\Omega_2$ (since it cannot be infinite). The same holds on
$\Omega_1$. Thus we can apply the Schwarz reflection to $\frac{Q}{\rmd
z^2}$ (which is up to now defined on $\Omega$), and we obtain a
holomorphic $2$-form $Q$ defined on $U\setminus\{0\}$.

In the same way, assuming that $D_2$ is the $x_1$-axis, we have
$\frac{g''}{g'}(z)\in\bar\R$ for $z\in\Omega_2$, and so 
$(\frac{g''}{g'})'-\frac{1}{2}(\frac{g''}{g'})^2$ is real or infinite
on $\Omega_2$. The same holds on $\Omega_1$, and we obtain a
meromorphic $2$-form $\rmS_z g$ defined on $U\setminus\{0\}$ by Schwarz
reflection. 
\end{proof}

\begin{lemma} \label{parameters}
Assume that $x$ has an end bounded by $(D_1,D_2)$. Let
$(g,\omega)$ be the Weierstrass data of $x$, and $Q$ its Hopf
differential. If $x$ has a helicoidal end bounded by $(D_1,D_2)$, then
there exist $A\in\R^*$ and $\alpha\in\R\setminus\Z$ such that 
\begin{equation} \label{hopfschwarzian}
Q\sim i\frac{A\alpha}{2\pi}z^{-2}\rmd z^2,\quad
\rmS_z g\sim \frac{1-\alpha^2}{2}z^{-2}\rmd z^2
\end{equation}
when $z\to 0$.

The couple $(A,\alpha)$ is then defined uniquely up to a sign. Moreover, if
$\pi\alpha_0$ denotes the angle of $\vec{u}_1$ and $-\vec{u}_2$ with
$\alpha_0\in(0,1)$, then we have either $\alpha\in\alpha_0+2\Z$ and
$A=-\rmD(D_1,D_2)$, or $-\alpha\in\alpha_0+2\Z$ and $A=\rmD(D_1,D_2)$.

We say that $x$ has a helicoidal end of parameters $(A,\alpha)$, and
that $\pi\alpha$ is the angle of the helicoidal end.
\end{lemma}

\begin{proof}
Since a direct isometry of $\R^3$ does not change $Q$ and $\rmS_z g$, we
can assume that 
$\vec{u}_1=(\cos(\pi\alpha_0),\sin(\pi\alpha_0),0)$,
$\vec{u}_2=-\vec{e}_1$, $D_1$ is the line $(0,0,-A)+\R\vec{u}_1$ and
$D_2$ the $x_1$-axis. Then we have $\vec{v}=-\vec{e}_3$, and so
$g(0)=0$ or $g(0)=\infty$ by lemma \ref{gaussend}. Moreover we have
$\rmD(D_1,D_2)=-A$. 

Set $h(z)=z^{-\alpha_0}g(z)$ for $z\in\Omega$. Then $h(z)\in i\bar{\R}$ if
$z\in\Omega_1$ (since $N(z)\in D_1^\perp$) or $z\in\Omega_2$ (since
$N(z)\in D_2^\perp$). 
Thus $h$ extends to a meromorphic map on
$U\setminus\{0\}$ by Schwarz reflection principle. Since $g$ has a
limit at $0$, $h$ has no essential singularity at $0$.

Thus there exist a integer $j$ and a nonzero real number $\rho$ such
that $g(z)\sim i\rho z^{\alpha_0+j}$. We set $\alpha=\alpha_0+j$. We
compute that $\rmS_z g\sim\frac{1-\alpha^2}{2}z^{-2}\rmd z^2$.

Let $h_1(z)=X_3(z)-i\frac{A}{\pi}\ln{z}$ for $z\in\Omega$, with $X_3$
as in section \ref{weierstrass}. Then, since $x_3=\re{X_3}$, we have
$h_1(z)\in i\R$ if $z\in\Omega_1$ or $z\in\Omega_2$. Thus $h_1$ extends to
a meromorphic map on $U\setminus\{0\}$ by Schwarz reflection
principle. Moreover, $x_3(z)=-<x(z),\vec{v}>$ is bounded in the
neighbourhood of $0$ by condition 4 in definition \ref{helicoidalend},
so $h_1$ is holomorphic at $0$.

Hence we have $\rmd X_3\sim i\frac{A}{\pi}z^{-1}\rmd z$, and so 
$Q=\frac{1}{2}\rmd X_3\frac{\rmd g}{g}\sim
i\frac{A\alpha}{2\pi}z^{-2}\rmd z^2$.

We now prove that the integer $j=\alpha-\alpha_0$ is even. We have 
$$\rmd(x_1+ix_2)=\bar\omega-g^2\omega
=\frac{\overline{\rmd X_3}}{2\bar g}-\frac{g\rmd X_3}{2}
\sim\frac{A}{2\pi}(\rho^{-1}\bar z^{-1-\alpha}\overline{\rmd z}
+\rho z^{-1+\alpha}\rmd z).$$ We set $z=t+i\tau$ with $t$ and $\tau$ real.
We have $<x,\vec{u}_1>=\re((x_1+ix_2)e^{-i\pi\alpha_0})$, so for $t<0$
we have
$$\frac{\partial}{\partial t}<x,\vec{u}_1>
\sim-\frac{A}{2\pi}\cos(\pi(\alpha-\alpha_0))
(\rho^{-1}|t|^{-1-\alpha}+\rho|t|^{-1+\alpha}).$$
Thus condition 1 in definition \ref{helicoidalend} implies that
$-A\rho\cos(\pi(\alpha-\alpha_0))>0$.
And we have $<x,-\vec{u}_2>=x_1=\re(x_1+ix_2)$, so for $t>0$ we have
$$\frac{\partial}{\partial t}<x,-\vec{u}_2>
\sim\frac{A}{2\pi}(\rho^{-1}t^{-1-\alpha}+\rho t^{-1+\alpha}).$$
Thus condition 2 in definition \ref{helicoidalend} implies that
$A\rho<0$. We conclude that 
$\cos(\pi(\alpha-\alpha_0))>0$, that is $\alpha-\alpha_0\in2\Z$.

Finally, it is clear that (\ref{hopfschwarzian}) defines $(A,\alpha)$
uniquely up to a sign.
\end{proof}

Let $x$ be an immersion having a helicoidal end of parameters
$(A,\alpha)$ bounded by two lines $D_1$ and $D_2$ that are as in the
proof of this lemma. Without loss of generality we can assume that
$\alpha\in\alpha_0+2\Z$, and thus $A=-\rmD(D_1,D_2)$. Then $\alpha>0$
means that the Gauss map at $0$ points down and that we turn in the
clockwise direction when we go from $D_1$ to $D_2$ on the minimal
surface, and $\alpha<0$ means that the Gauss map at $0$ points up and
that we turn in the counter-clockwise direction when we go from $D_1$
to $D_2$ on the minimal surface. On the other hand, $A>0$ means that
$D_1$ lies below $D_2$, and $A<0$ means that $D_1$ lies above
$D_2$. Thus, $A\alpha>0$ means that we go down when we turn in the
counter-clockwise direction on the minimal surface, and $A\alpha<0$
means that we go up when we turn in the counter-clockwise direction on
the minimal surface. This last fact remains true if
$-\alpha\in\alpha_0+2\Z$. Hence we say that $x$ has a left-helicoidal
end (respectively a right-helicoidal end) if $A\alpha>0$ (respectively
$A\alpha<0$) (see figures \ref{pospos}, \ref{posneg}, \ref{negpos} and
\ref{negneg}).

\begin{figure}[htbp]
\begin{center}
\input{pospos.pstex_t}
\caption{$A>0$ and $\alpha>0$ (left-helicoidal end).}
\label{pospos}
\end{center}
\end{figure}
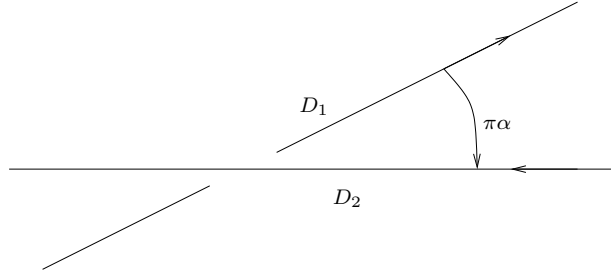

\begin{figure}[htbp]
\begin{center}
\input{posneg.pstex_t}
\caption{$A<0$ and $\alpha>0$ (right-helicoidal end).}
\label{posneg}
\end{center}
\end{figure}

\begin{figure}[htbp]
\begin{center}
\input{negpos.pstex_t}
\caption{$A>0$ and $\alpha<0$ (right-helicoidal end).}
\label{negpos}
\end{center}
\end{figure}

\begin{figure}[htbp]
\begin{center}
\input{negneg.pstex_t}
\caption{$A<0$ and $\alpha<0$ (left-helicoidal end).}
\label{negneg}
\end{center}
\end{figure}

\begin{rem}
This definition and these lemmas extend for ends at a point
$z_1\in\R$. 

They also extend for end at $\infty$ using the change of
parameters $\zeta=-z^{-1}$, which maps $\{\im{z}>0\}$ onto itself. We
get
$$Q\sim i\frac{A\alpha}{2\pi}z^{-2}\rmd z^2,\quad
\rmS_z g\sim \frac{1-\alpha^2}{2}z^{-2}\rmd z^2$$ when $z\to\infty$
(because $\rmS_zg=\rmS_\zeta g+\rmS_z\zeta$ and $\rmS_z\zeta=0$).
\end{rem}

\section[Minimal surfaces bounded by three lines]{Minimal surfaces
bounded by three lines with helicoidal ends} 
\label{minimalsurface}

In this section we will study minimal disks
bounded by three lines with three helicoidal ends when the triple of lines
belong to $\cD$, which is a generic property.

\subsection{Geometric configuration} \label{geometricconfiguration}

An element of $\cD$ has a representant that is described as follows.

Let $D_1$ be the horizontal line oriented by the vector
$\vec{u}_1=(\cos(\pi\alpha_0),\sin(\pi\alpha_0),0)$ for some
$\alpha_0\in(0,1)$, 
$D_2$ the $x_1$-axis oriented by the vector $\vec{u}_2=-\vec{e}_1$, 
and $D_3$ the line oriented by the vector $\vec{u}_3=
(\cos(\pi\gamma')\sin\kappa,-\sin(\pi\gamma')\sin\kappa,\cos\kappa)$
for some $\gamma'\in\R$ and $\kappa\in\R$ (the number $\pi\gamma'$ is
the angle of the projections of $D_2$ and $D_3$ on the horizontal
plane, except if $D_3$ is vertical, in what case it can take any value).

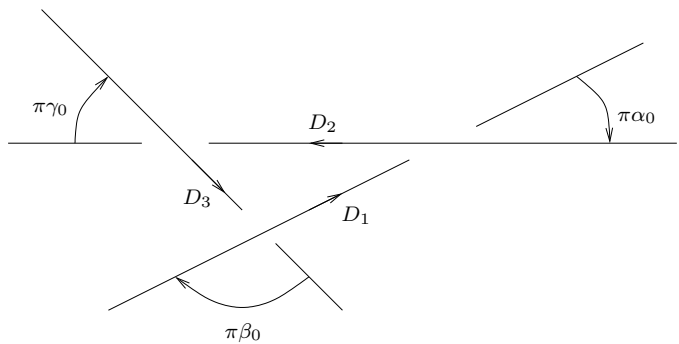
\begin{figure}[htbp]
\begin{center}
\input{3lines.pstex_t}
\caption{Three lines in generic position.}
\label{configuration}
\end{center}
\end{figure}

The number $\pi\alpha_0$ is the geometric angle of $\vec{u}_1$ and
$-\vec{u}_2$. Let us denote by $\pi\beta_0$ with
$\beta_0\in(0,1)$ the geometric angle of $\vec{u}_3$ and
$-\vec{u}_1$, and by $\pi\gamma_0$ with
$\gamma_0\in(0,1)$ the geometric angle of $\vec{u}_2$ and
$-\vec{u}_3$ (see figure \ref{configuration}).

We denote by $\vec{v}_0$, $\vec{v}_1$ and $\vec{v}_{\infty}$ the
vectors associated to $(D_1,D_2)$, $(D_2,D_3)$ and $(D_3,D_1)$ (see
definition \ref{associatedvector}):
$$\vec{v}_0=-\frac{\vec{u}_1\times\vec{u}_2}{||\vec{u}_1\times\vec{u}_2||}
=-\vec{e}_3,\quad
\vec{v}_1=-\frac{\vec{u}_2\times\vec{u}_3}{||\vec{u}_2\times\vec{u}_3||},
\quad
\vec{v}_{\infty}
=-\frac{\vec{u}_3\times\vec{u}_1}{||\vec{u}_3\times\vec{u}_1||}.$$

We set $A=-\rmD(D_1,D_2)$, $B=-\rmD(D_3,D_1)$ and $C=-\rmD(D_2,D_3)$.
Finally we denote by $\varepsilon_0$ the sign of
$\det(\vec{u}_1,\vec{u}_2,\vec{u}_3)$.

\begin{prop} \label{bijection}
The map
$$L:(D_1,D_2,D_3)\mapsto(\alpha_0,\gamma_0,\beta_0,-A,-C,-B,\varepsilon_0)$$
is a bijection from $\cD$ onto 
$\cK\times\R^*\times\R^*\times\R^*\times\{1,-1\}$ where
$\cK$ is the set of the triples $(\alpha_0,\gamma_0,\beta_0)\in\R^3$
satisfying
\begin{equation} 
\label{sphericalangles}
\begin{array}{cc}
\alpha_0+\beta_0+\gamma_0>1, & -\alpha_0+\beta_0+\gamma_0<1,\\
\alpha_0-\beta_0+\gamma_0<1, & \alpha_0+\beta_0-\gamma_0<1.
\end{array}
\end{equation}
\end{prop}

\begin{proof}
The fact that $(\alpha_0,\gamma_0,\beta_0)\in\cK$ is a consequence of
Gauss-Bonnet formula applied to the spherical triangles on $\s^2$
bounded by the circles $D_1^\perp$, $D_2^\perp$ and $D_3^\perp$. 

Conversely, let $(\alpha_0,\gamma_0,\beta_0)\in\cK$, $A,B,C\in\R^*$
and $\varepsilon_0\in\{1,-1\}$. Then there exists a spherical triangle
of angles $\pi\alpha_0$, $\pi\beta_0$ and $\pi\gamma_0$. The three
corresponding oriented circles define unit vectors $\vec{u}_1$,
$\vec{u}_2$ and $\vec{u}_3$ uniquely up to a direct isometry of
$\R^3$. If the sign of $\det(\vec{u}_1,\vec{u}_2,\vec{u}_3)$ 
is not equal to $\varepsilon_0$, then we replace these vectors by
their images by an indirect isometry of $\R^3$ (which does not change
the angles of the spherical triangle). Up to now, $\vec{u}_1$,
$\vec{u}_2$ and $\vec{u}_3$ are uniquely determined.

Now we consider three lines $D_1$, $D_2$ and $D_3$ in $\R^3$ oriented
by  $\vec{u}_1$, $\vec{u}_2$ and $\vec{u}_3$. We translate $D_2$ in
the direction of $\vec{u}_1\times\vec{u}_2$ in order that
$\rmD(D_1,D_2)=-A$. Then we translate $D_3$ in the direction of
$\vec{u}_1\times\vec{u}_3$ in order that $\rmD(D_3,D_1)=-B$. Finally
we translate $D_3$ in the direction of $\vec{u}_1$ in order that
$\rmD(D_2,D_3)=-C$ (this operation does not change
$\rmD(D_3,D_1)$). The lines $D_1$, $D_2$ and $D_3$ are determined
uniquely up to a direct isometry of $\R^3$. This completes the proof.
\end{proof}

A precise study of the space of triples of lines in generic position,
and in particular a more detailed proof of proposition
\ref{bijection}, can be found in \cite{balser}.

\begin{defn} \label{dual}
Two triples in $\cD$ are called dual configurations if their parameters
only differ by the sign of $\varepsilon_0$.
\end{defn}

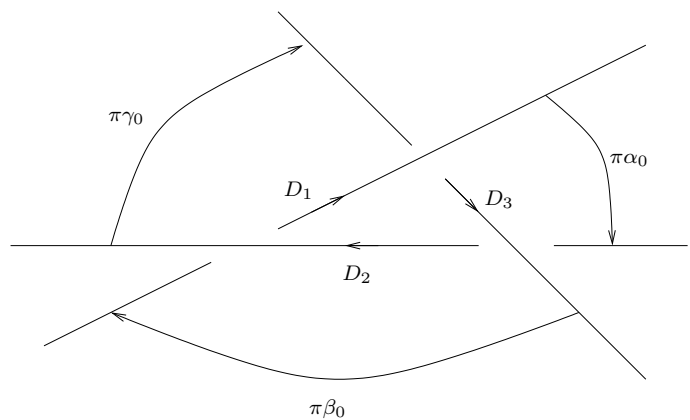
\begin{figure}[htbp]
\begin{center}
\input{dual.pstex_t}
\caption{The dual configuration of that of figure \ref{configuration}.}
\label{figuredual}
\end{center}
\end{figure}

The dual configuration of that of figure \ref{configuration} is shown
on figure \ref{figuredual}; in both configurations the lines $D_1$
and $D_2$ are horizontal, but $D_3$ ``goes down'' on figure
\ref{configuration} and ``goes up'' on figure \ref{figuredual}.

\begin{rem}
An indirect isometry of $\R^3$ changes $A$, $B$, $C$ and
$\varepsilon_0$ into their opposites. The dual configuration of a
triple is not the image of this triple by a symmetry, but the
directions of the straight lines are symmetric.
\end{rem}

We set $\beta'=1-\alpha-\gamma'$ (the number $\pi\beta'$ is
the angle of the projections of $D_1$ and $D_3$ on the horizontal
plane, except if $D_3$ is vertical, in what case it can take any value).

Since $\vec{v}_1$ and $-\vec{e}_3$ are normal to the $x_1$-axis, there
exists a rotation $R$ about the oriented $x_1$-axis that maps 
$\vec{v}_1$ onto $-\vec{e}_3$; we denote by $\theta$ its angle. We
have $R(\vec{u}_2)=\vec{u}_2$ and
$R(\vec{u}_3)=(\cos(\pi\gamma_0),-\sin(\pi\gamma_0),0)$.
In the same way, $\vec{v}_{\infty}$ and $-\vec{e}_3$ are normal to
$D_1$, so there exists a rotation $\hat R$ about the oriented line
$D_1$ that maps $\vec{v}_{\infty}$ onto $-\vec{e}_3$; we denote by
$\hat\theta$ its angle. We denote by $T$ the rotation of angle
$-\pi\alpha_0$ with respect to the $x_3$-axis. We have
$T\circ\hat R(\vec{u}_1)=\vec{e}_1$,
$T\circ\hat R(\vec{u}_3)=(\cos(\pi\beta_0),-\sin(\pi\beta_0),0)$ and
$T\circ\hat R(\vec{v}_\infty)=-\vec{e}_3$.

Finally we set $t=\tan\frac{\theta}{2}$ and
$\hat{t}=\tan\frac{\hat\theta}{2}$.

\begin{lemma} \label{costheta}
We have $$\cos\theta=\frac{\cos(\pi\beta_0)+\cos(\pi\alpha_0)\cos(\pi\gamma_0)}
{\sin(\pi\alpha_0)\sin(\pi\gamma_0)},$$
$$\cos\hat\theta=\frac{\cos(\pi\gamma_0)+\cos(\pi\alpha_0)\cos(\pi\beta_0)}
{\sin(\pi\alpha_0)\sin(\pi\beta_0)}.$$
\end{lemma}

\begin{proof}
We notice that the numbers $\sin(\pi\alpha_0)$, $\sin(\pi\beta_0)$ and
$\sin(\pi\gamma_0)$ are positive.

We have $$\cos(\pi\gamma_0)=<\vec{u}_2,-\vec{u}_3>
=\cos(\pi\gamma')\sin\kappa,$$
$$\sin(\pi\gamma_0)=\sqrt{1-\cos^2(\pi\gamma_0)}=
\sqrt{\sin^2(\pi\gamma')\sin^2\kappa+\cos^2\kappa},$$
\begin{eqnarray*}
\cos(\pi\beta_0)& = & <\vec{u}_3,-\vec{u}_1> \\
& = & -\cos(\pi\alpha_0)\cos(\pi\gamma')\sin\kappa
+\sin(\pi\alpha_0)\sin(\pi\gamma')\sin\kappa.
\end{eqnarray*}
We compute that
$$\vec{u}_1\times\vec{u}_2=\sin(\pi\alpha_0)\vec{e}_3,\quad 
\vec{u}_2\times\vec{u}_3=(0,\cos\kappa,\sin(\pi\gamma')\sin\kappa).$$
Thus we have
$$\cos\theta=\frac{<\vec{u}_1\times\vec{u}_2,\vec{u}_2\times\vec{u}_3>}
{||\vec{u}_1\times\vec{u}_2||\cdot||\vec{u}_2\times\vec{u}_3||}
=\frac{\sin(\pi\gamma')\sin\kappa}
{\sqrt{\sin^2(\pi\gamma')\sin^2\kappa+\cos^2\kappa}}.$$
Finally we get
$$\cos(\pi\beta_0)=
-\cos(\pi\alpha_0)\cos(\pi\gamma_0)
+\sin(\pi\alpha_0)\sin(\pi\gamma_0)\cos\theta.$$
This proves the first formula.

And we have
$$\cos(\pi\beta_0)=\cos(\pi\beta')\sin\kappa,$$
$$\sin(\pi\beta_0)=\sqrt{1-\cos^2(\pi\beta_0)}=
\sqrt{\sin^2(\pi\beta')\sin^2\kappa+\cos^2\kappa}.$$
We compute that $$\vec{u}_3\times\vec{u}_1=(-\sin(\pi\alpha_0)\cos\kappa,
\cos(\pi\alpha_0)\cos\kappa,\sin(\pi\beta')\sin\kappa).$$
Thus we have
$$\cos\hat\theta=\frac{<\vec{u}_1\times\vec{u}_2,\vec{u}_3\times\vec{u}_1>}
{||\vec{u}_1\times\vec{u}_2||\cdot||\vec{u}_3\times\vec{u}_1||}
=\frac{\sin(\pi\beta')\sin\kappa}
{\sqrt{\sin^2(\pi\beta')\sin^2\kappa+\cos^2\kappa}}.$$
Finally we have
\begin{eqnarray*}
\cos(\pi\gamma_0) & = & -\cos(\pi(\alpha_0+\beta'))\sin\kappa \\
& = & -\cos(\pi\alpha_0)\cos(\pi\beta')\sin\kappa
+\sin(\pi\alpha_0)\sin(\pi\beta')\sin\kappa \\
& = & -\cos(\pi\alpha_0)\cos(\pi\beta_0)
+\sin(\pi\alpha_0)\sin(\pi\beta_0)\cos\hat\theta.
\end{eqnarray*}
This proves the second formula.
\end{proof}

\begin{lemma} \label{signs}
The signs of $\cos\kappa$, $\sin\theta$, $\sin\hat\theta$, $t$ and
$\hat t$ are equal to $\varepsilon_0$.
\end{lemma}

\begin{proof}
Since $\sin\theta=\frac{2t}{1+t^2}$, $t$ and $\sin\theta$ have the
same sign. In the same way $\hat t$ and $\sin\hat\theta$ have the same sign. 

By definition of $\theta$ we have
$\vec{v}_1\times(-\vec{e}_3)=\sin\theta\vec{e}_1$. We compute that
$$\vec{v}_1\times(-\vec{e}_3)=
\frac{\cos\kappa}{\sqrt{\sin^2(\pi\gamma')\sin^2\kappa+\cos^2\kappa}}
\vec{e}_1,$$
so $\cos\kappa$ and $\sin\theta$ have the same sign.

In the same way we have $\vec{v}_{\infty}\times(-\vec{e}_3)
=\sin\hat\theta\vec{u}_1$, so
$$\sin\hat\theta=\det(\vec{v}_{\infty},-\vec{e}_3,\vec{u}_1)
=\frac{\cos\kappa}{\sqrt{\sin^2(\pi\beta')\sin^2\kappa+\cos^2\kappa}}.$$

Finally we have
$\det(\vec{u}_1,\vec{u}_2,\vec{u}_3)=\sin(\pi\alpha_0)\cos\kappa$.
\end{proof}

\subsection{The Hopf differential and the spinor data}
\label{hopfandspinors}

Proceeding as for lemma \ref{hopfschwarziansymmetry}, we get the
following result.

\begin{lemma} \label{hopfschwarzianmeromorphic}
Let $x:\Sigma\to\R^3$ be a conformal minimal immersion bounded by
three straight lines. Let $(g,\omega)$ be its Weierstrass data. Then
its Hopf differential $Q$ extends to a holomorphic $2$-form on
$\C\setminus\{0,1\}$ and the Schwarzian derivative $\rmS_z g$ of its
Gauss map extends to a meromorphic $2$-form on $\C\setminus\{0,1\}$.
\end{lemma}

From now on we consider a triple of lines $(D_1,D_2,D_3)\in\cD$. Let
$$L(D_1,D_2,D_3)=(\alpha_0,\gamma_0,\beta_0,-A,-C,-B,\varepsilon_0),$$
$\alpha\in\alpha_0+2\Z$, $\beta\in\beta_0+2\Z$ and
$\gamma\in\gamma_0+2\Z$. We assume that $x:\Sigma\to\R^3$ is a
conformal minimal immersion bounded by $(D_1,D_2,D_3)$ and having
helicoidal ends of parameters $(A,\alpha)$, $(B,\beta)$ and
$(C,\gamma)$ at $0$, $\infty$ and $1$ respectively. We denote by
$(g,\omega)$ its Weierstrass data.

\begin{prop}[Riemann, \cite{riemann}] \label{hopf}
Then the Hopf differential of $x$ is
\begin{equation} \label{eqQ}
Q=iz^{-2}(z-1)^{-2}\varphi(z)\rmd z^2
\end{equation}
where
\begin{equation} \label{eqphi}
\varphi(z)=\frac{B\beta}{2\pi}z(z-1)
-\frac{A\alpha}{2\pi}(z-1)+\frac{C\gamma}{2\pi}z.
\end{equation}
\end{prop}

\begin{proof}
At $z=0$ we have $Q\sim i\frac{A\alpha}{2\pi}z^{-2}\rmd z^2$.
At $z=1$ we have $Q\sim i\frac{C\gamma}{2\pi}(z-1)^{-2}\rmd z^2$.
At $z=\infty$ we have $Q\sim i\frac{B\beta}{2\pi}z^{-2}\rmd z^2$.
Hence the map $\varphi(z)=z^2(z-1)^2\frac{Q}{i\rmd z^2}$ has no
singularity on $\C$, and we have $\varphi(z)=\rmO(z^2)$ when
$z\to\infty$, so $\varphi$ is a polynomial of degree less than or
equal to $2$.
Finally we compute that $\varphi$ has the announced expression.
\end{proof}

\begin{lemma} \label{lemmazerophi}
The polynomial $\varphi$ defined by (\ref{eqphi}) has two nonreal
conjugate roots if and only
if $A\alpha$, $B\beta$ and $C\gamma$ have the same sign,
$\sqrt{|A\alpha|}<\sqrt{|B\beta|}+\sqrt{|C\gamma|}$, 
$\sqrt{|B\beta|}<\sqrt{|A\alpha|}+\sqrt{|C\gamma|}$ and
$\sqrt{|C\gamma|}<\sqrt{|A\alpha|}+\sqrt{|B\beta|}$.

It has a double real root if and only if $A\alpha$, $B\beta$ and
$C\gamma$ have the same sign, and
$\sqrt{|A\alpha|}=\sqrt{|B\beta|}+\sqrt{|C\gamma|}$, 
$\sqrt{|B\beta|}=\sqrt{|A\alpha|}+\sqrt{|C\gamma|}$ or
$\sqrt{|C\gamma|}=\sqrt{|A\alpha|}+\sqrt{|B\beta|}$.

It has two distinct real roots in all other cases.
\end{lemma}

\begin{proof}
The discriminant of $\varphi$ is $\frac{\delta}{4\pi^2}$ where
$$\delta=A^2\alpha^2+B^2\beta^2+C^2\gamma^2-2AB\alpha\beta
-2AC\alpha\gamma-2BC\beta\gamma.$$ Thus the three cases in the lemma
correspond respectively to $\delta<0$, $\delta=0$ and $\delta>0$.
The expression $\delta$ is a polynomial in the variable $C$, whose
discriminant is equal to $16AB\alpha\beta\gamma^2$. If $A\alpha
B\beta<0$, then $\delta(C)>0$ for all $C\in\R^*$. 

We assume that $A\alpha>0$ and $B\beta>0$. Then we have $\delta(C)<0$
if and only if $$(\sqrt{A\alpha}-\sqrt{B\beta})^2<C\gamma
<(\sqrt{A\alpha}+\sqrt{B\beta})^2.$$ This condition is not satisfied if
$C\gamma<0$, and if $C\gamma>0$ it is satisfied if and only if
$\sqrt{C\gamma}>\sqrt{A\alpha}-\sqrt{B\beta}$, 
$\sqrt{C\gamma}>\sqrt{B\beta}-\sqrt{A\alpha}$ and 
$\sqrt{C\gamma}<\sqrt{A\alpha}+\sqrt{B\beta}$. And we have
$\delta(C)=0$ if and only if $C>0$, and
$\sqrt{C\gamma}=\sqrt{A\alpha}-\sqrt{B\beta}$,  
$\sqrt{C\gamma}=\sqrt{B\beta}-\sqrt{A\alpha}$ or 
$\sqrt{C\gamma}=\sqrt{A\alpha}+\sqrt{B\beta}$.

We deal with the case where $A\alpha<0$ and $B\beta<0$ in the same way.
\end{proof}

We set $$\zeta(z)=\int_0^z\varphi(\tau)\rmd\tau.$$ This is a local
diffeomorphism except at the zeros of $\varphi$. The Schwarzian
derivative of $g$ in the ``coordinate'' $\zeta$ satisfies
$$\rmS_\zeta g=\rmS_zg-\rmS_z\zeta.$$ We set
\begin{eqnarray*}
\Theta=\frac{\rmS_\zeta g}{\rmd z^2}
& = & \varphi\left(\frac{1}{g'}\left(\frac{g'}{\varphi}\right)'\right)'
-\frac{\varphi^2}{2}
\left(\frac{1}{g'}\left(\frac{g'}{\varphi}\right)'\right)^2 \\
& = & \varphi^2\left(\left(\frac{g''}{g'}\right)'
-\frac{1}{2}\left(\frac{g''}{g'}\right)^2
-\left(\frac{\varphi''}{\varphi}\right)'
+\frac{1}{2}\left(\frac{\varphi'}{\varphi}\right)^2\right).
\end{eqnarray*}

Let $(\xi_1,\xi_2)$ be the spinor data of $x$. We
set $$\xi_1=z^{-1}(z-1)^{-1}k_1\sqrt{\rmd z},\quad
\xi_2=z^{-1}(z-1)^{-1}k_2\sqrt{\rmd z}.$$ Then $k_1$ and $k_2$ are
holomorphic functions on $\Sigma$, and we have
\begin{equation} \label{defk1k2}
\begin{array}{c}
g=\frac{k_2}{k_1},\quad\omega=z^{-2}(z-1)^{-2}k_1^2\rmd z, \\
Q=z^{-2}(z-1)^{-2}(k_1k_2'-k_1'k_2)\rmd z^2.
\end{array}
\end{equation}
We will abusively call $k_1$ and $k_2$ the spinors associated to
$(g,\omega)$. 

\begin{lemma}[Riemann, \cite{riemann}] \label{lemmak1k2}
The functions $k_1$ and $k_2$ satisfy the following relation on $\Sigma$:
\begin{equation} \label{eqk1k2}
k_1k_2'-k_1'k_2=i\varphi.
\end{equation}
They are solutions on $\Sigma$ of the following
differential equation:
\begin{equation} \label{eqk}
k''-\frac{\varphi'}{\varphi}k'+\frac{\Theta}{2}k=0.
\end{equation}
\end{lemma}

\begin{proof}
Equality (\ref{eqk1k2}) follows from (\ref{eqQ}).

Since $k_1^2=\frac{i\varphi}{g'}$, we have
$$2\frac{k_1'}{k_1}=\frac{\varphi'}{\varphi}-\frac{g''}{g'}.$$ On the
other hand we have
$$\frac{1}{g'}\left(\frac{g'}{\varphi}\right)'
=\frac{1}{g'}\left(\frac{g''}{\varphi}-\frac{\varphi'g'}{\varphi^2}\right)
=-2\frac{k_1'}{\varphi k_1}.$$
So $$\left(\frac{1}{g'}\left(\frac{g'}{\varphi}\right)'\right)'
=-2\frac{k_1''}{\varphi k_1}+2\frac{{k_1'}^2}{\varphi k_1^2}
+2\frac{\varphi'k_1'}{\varphi^2k_1},$$ and
$$\Theta=\varphi
\left(\frac{1}{g'}\left(\frac{g'}{\varphi}\right)'\right)'
-\frac{\varphi^2}{2}\left(\frac{1}{g'}\left(\frac{g'}{\varphi}\right)'\right)^2
=-2\frac{k_1''}{k_1}+2\frac{\varphi'k_1'}{\varphi k_1}.$$
Thus $k_1$ is solution of equation (\ref{eqk}).

We also have $$2\frac{k_2'}{k_2}=
2\frac{g'}{g}+\frac{\varphi'}{\varphi}-\frac{g''}{g'}.$$
So $$\frac{1}{g'}\left(\frac{g'}{\varphi}\right)'
=2\frac{g'}{\varphi g}-2\frac{k_2'}{\varphi k_2},$$
$$\left(\frac{1}{g'}\left(\frac{g'}{\varphi}\right)'\right)'
=2\frac{g''}{\varphi g}-2\frac{g'^2}{\varphi g^2}
-2\frac{\varphi'g'}{\varphi^2g}-2\frac{k_2''}{\varphi k_2}
+2\frac{{k_2'}^2}{\varphi k_2^2}+2\frac{\varphi'k_2'}{\varphi^2k_2},$$
$$\Theta=2\frac{g''}{g}-4\frac{g'^2}{g^2}
-2\frac{\varphi'g'}{\varphi g}-2\frac{k_2''}{k_2}
+2\frac{\varphi'k_2'}{\varphi k_2}+4\frac{g'k_2'}{gk_2}
=-2\frac{k_2''}{k_2}+2\frac{\varphi'k_2'}{\varphi k_2}.$$
Thus $k_2$ is solution of equation (\ref{eqk}).
\end{proof}

\begin{lemma} \label{Theta}
We have $$\Theta(z)=\frac{\Phi(z)}{z^2(z-1)^2}
+\frac{\Lambda(z)}{z(z-1)\varphi(z)}+\frac{2\varphi''}{\varphi(z)}$$ where
\begin{equation} \label{Phi}
\Phi(z)=\frac{1-\beta^2}{2}z(z-1)
-\frac{1-\alpha^2}{2}(z-1)+\frac{1-\gamma^2}{2}z
\end{equation}
and where $\Lambda$ is an affine function.
\end{lemma}

\begin{proof}
The form $\rmS_z\zeta$ is meromorphic on $\bar\C$, with double poles
at $\infty$ and the roots of $\varphi$.
Since $\Theta\rmd z^2=\rmS_zg-\rmS_z\zeta$, by section
\ref{helicoidalends} and lemma \ref{hopfschwarzianmeromorphic} the
function $\Theta$ is meromorphic on $\bar\C$, its possible poles are
$0$, $1$, $\infty$ and the zeros of $\varphi$, and these poles are at
most double.

By lemma \ref{parameters} we have
$\rmS_zg\sim\frac{1-\alpha^2}{2}z^{-2}\rmd z^2$ when $z\to 0$. On the
other hand, $\rmS_z\zeta$ is holomorphic at $0$ since $0$ is not a
zero of $\varphi$. Thus we have
$\Theta(z)\sim\frac{1-\alpha^2}{2}z^{-2}$ when $z\to 0$. In
the same way we have
$\Theta(z)\sim\frac{1-\gamma^2}{2}(z-1)^{-2}$ when $z\to 1$ and 
$\Theta(z)\sim\frac{9-\beta^2}{2}z^{-2}\rmd z^2$ when
$z\to\infty$ (since $\rmS_z\zeta\sim-4z^{-2}\rmd z^2$). At a root of
$\varphi$, since $g$ and $\varphi$ have 
branch points of the same order, the order $-2$ terms in $\rmS_zg$
and $\rmS_z\zeta$ at this root are equal, and so the order of $\Theta$
is greater than or equal to $-1$.

We have $\frac{2\varphi''}{\varphi(z)}\sim 4z^{-2}$ when $z\to\infty$
(and $\varphi''$ is a constant). Consequently the function
$\Lambda=z(z-1)\varphi(\Theta-z^{-2}(z-1)^{-2}\Phi)
-2z(z-1)\varphi''$ is holomorphic 
on $\C$, and we have $\Lambda(z)=\rmO(z)$ when $z\to\infty$, so it is
an affine function.
\end{proof}

\begin{lemma} \label{fourschwarzian}
If the roots of $\varphi$ are distinct, then there are at most four
possibilities for the function $\Lambda$. 
\end{lemma}

\begin{proof}
Equation (\ref{eqk}) has regular singularities at $0$, $1$, $\infty$ and
the roots $a_1$ and $a_2$ of $\varphi$ (see \cite{ww}, paragraph
10.3). Moreover, at least one of the roots of $\varphi$ lies in
$\Sigma$, for example $a_1$.

Since $a_1\neq a_2$, the exponents of
equation (\ref{eqk}) at $a_1$ are $0$ and $2$, and since $k_1$ and
$k_2$ are well-defined on $\Sigma$, the solutions of (\ref{eqk}) have
no logarithmic term at $a_1$. Thus, in the neighbourhood of $a_1$,
equation (\ref{eqk}) has a solution having the following form:
$$\sum_{n=0}^\infty\lambda_n(z-a_1)^n$$ with $\lambda_0\neq 0$.
Writing
$\Theta(z)=\psi_{-1}(z-a_1)^{-1}+\psi_0+\rmO(z-a_1)$, since
$\frac{\varphi'(z)}{\varphi(z)}=\frac{1}{z-a_1}+\frac{1}{a_1-a_2}
+\rmO(z-a_1)$, reporting in equation (\ref{eqk}) we get 
$$-\lambda_1+\frac{\psi_{-1}}{2}\lambda_0=0,\quad
\left(\frac{\psi_{-1}}{2}-\frac{1}{a_1-a_2}\right)\lambda_1
+\frac{\psi_0}{2}\lambda_0=0.$$
Hence we get 
\begin{equation} \label{psi-1psi0}
\psi_{-1}\left(\frac{\psi_{-1}}{2}-\frac{1}{a_1-a_2}\right)+\psi_0=0.
\end{equation}

We also compute using (\ref{psi-1psi0}) that 
$$\Lambda(z)=\frac{B\beta}{2\pi}(-4a_1(a_1-1)+m_1\psi_{-1})
+\frac{B\beta}{2\pi}(m_2+m_3\psi_{-1}+m_4\psi_{-1}^2)(z-a_1)$$
with
$$m_1=a_1(a_1-1)(a_1-a_2),\quad
m_2=-\frac{(a_1-a_2)\Phi(a_1)}{a_1(a_1-1)}-4(2a_1-1),$$
$$m_3=2a_1(a_1-1)+(2a_1-1)(a_1-a_2),\quad
m_4=-\frac{m_1}{2}.$$

Assume that $a_1$ and $a_2$ are distinct real roots. Then
$a_2\in\Sigma$ and we can apply the same argument at $a_2$ with a coefficient
$\tilde\psi_{-1}$ analogous to $\psi_{-1}$: we have 
$$\Lambda(z)=\frac{B\beta}{2\pi}(-4a_2(a_2-1)+\tilde m_1\tilde\psi_{-1})
+\frac{B\beta}{2\pi}(\tilde m_2+\tilde m_3\tilde\psi_{-1}
+\tilde m_4\tilde\psi_{-1}^2)(z-a_2)$$
where $\tilde m_j$ has the same expression as $m_j$ exchanging $a_1$
and $a_2$. Identifying these two expressions of $\Lambda$, we get
$\tilde m_2+\tilde m_3\tilde\psi_{-1}+\tilde m_4\tilde\psi_{-1}^2=
m_2+m_3\psi_{-1}+m_4\psi_{-1}^2$ and
$m_1\psi_{-1}-\tilde m_1\tilde\psi_{-1}
-4a_1(a_1-1)+4a_2(a_2-1)
=(a_1-a_2)(m_2+m_3\psi_{-1}+m_4\psi_{-1}^2)$. Setting
$R=a_1(a_1-1)\psi_{-1}$ and $\tilde R=a_2(a_2-1)\tilde\psi_{-1}$, we get
\begin{eqnarray} \label{rtilder}
\left\{ \begin{array}{ccc}
0 & = & -4-\frac{\Phi(a_1)}{a_1(a_1-1)}+
\left(\frac{1}{a_1-a_2}+\frac{1}{a_1}+\frac{1}{a_1-1}\right)R \\
& & -\frac{1}{a_1-a_2}\tilde R-\frac{1}{2a_1(a_1-1)}R^2 \\
0 & = & -4-\frac{\Phi(a_2)}{a_2(a_2-1)}+
\left(\frac{1}{a_2-a_1}+\frac{1}{a_2}+\frac{1}{a_2-1}\right)\tilde R \\
& & -\frac{1}{a_2-a_1}R-\frac{1}{2a_2(a_2-1)}\tilde R^2.
\end{array} \right.
\end{eqnarray}

Assume that $a_1$ and $a_2$ are complex conjugate roots. Since
$\Theta(z)\in\R$ when $z\in\R\setminus\{0,1\}$, $\Lambda$ must have
real coefficients, so we have
$\im(m_2+m_3\psi_{-1}+m_4\psi_{-1}^2)=0$ and
$\im(-4a_1(a_1-1)+m_1\psi_{-1}-a_1(m_2+m_3\psi_{-1}+m_4\psi_{-1}^2))=0$.
Setting $R=a_1(a_1-1)\psi_{-1}$ and $\tilde R=\bar R$, since
$a_2=\overline{a_1}$, this is equivalent to system (\ref{rtilder}).

System (\ref{rtilder}) has at most four solutions $(R,\tilde R)$. Thus
there are at most four possible functions $\Lambda$.
\end{proof}

\begin{lemma} \label{fourschwarziandoubleroot}
If the polynomial $\varphi$ has a double real root, then there are at
most three possibilities for the function $\Lambda$. 
\end{lemma}

\begin{proof}
Equation (\ref{eqk}) has regular singularities at $0$, $1$, $\infty$ and
the root $a_1$ of $\varphi$. Moreover, $a_1$ lies in $\Sigma$.
Since $a_1$ is a double root of $\varphi$, the exponents of equation
(\ref{eqk}) at $a_1$ are $0$ and $3$, and since $k_1$ and 
$k_2$ are well-defined on $\Sigma$, the solutions of (\ref{eqk}) have
no logarithmic term at $a_1$. Thus, in the neighbourhood of $a_1$,
equation (\ref{eqk}) has a solution having the following form:
$$\sum_{n=0}^\infty\lambda_n(z-a_1)^n$$ with $\lambda_0\neq 0$.
Writing
$$\Theta(z)=\psi_{-1}(z-a_1)^{-1}+\psi_0+\psi_1(z-a_1)+\rmO((z-a_1)^2),$$
since $\frac{\varphi'(z)}{\varphi(z)}=\frac{2}{z-a_1}$, reporting in
equation (\ref{eqk}) we get  
$$-2\lambda_1+\frac{\psi_{-1}}{2}\lambda_0=0,\quad
-2\lambda_2+\frac{\psi_{-1}}{2}\lambda_1+\frac{\psi_0}{2}\lambda_0=0,$$
$$\frac{\psi_{-1}}{2}\lambda_2+\frac{\psi_0}{2}\lambda_1
+\frac{\psi_1}{2}\lambda_0=0.$$
Hence we get 
\begin{equation} \label{psi-1psi0psi1}
\frac{\psi_{-1}^3}{16}+\frac{\psi_{-1}\psi_0}{2}+\psi_1=0.
\end{equation}

We also compute that
\begin{eqnarray*}
\Lambda(z)& = & -4\frac{B\beta}{2\pi}a_1(a_1-1)
+\frac{B\beta}{2\pi}n_1(z-a_1) \\
& & +\frac{B\beta}{2\pi}n_2(z-a_1)^2
+\frac{B\beta}{2\pi}n_3(z-a_1)^3+\rmO((z-a_1)^4)
\end{eqnarray*}
with $$n_1=a_1(a_1-1)\psi_{-1}-4(2a_1-1),$$
$$n_2=a_1(a_1-1)\left(\psi_0-\frac{\Phi(a_1)}{a_1^2(a_1-1)^2}\right)
+(2a_1-1)\psi_{-1}-4,$$
\begin{eqnarray*}
n_3 & = & a_1(a_1-1)\psi_1-\frac{\Phi'(a_1)}{a_1(a_1-1)}
+2\frac{\Phi(a_1)}{a_1(a_1-1)}
\left(\frac{1}{a_1}+\frac{1}{a_1-1}\right) \\
& & +(2a_1-1)\left(\psi_0-\frac{\Phi(a_1)}{a_1^2(a_1-1)^2}\right)
+\psi_{-1}.
\end{eqnarray*}
On the other hand, $\Lambda$ is an affine function, so we have
$n_2=n_3=0$, and using (\ref{psi-1psi0psi1}) we obtain that
$\psi_{-1}$ is solution of a degree $3$ polynomial equation.
Thus there are at most three possible functions $\Lambda$.
\end{proof}

All that has been done up to now does not depend on $\varepsilon_0$,
{\it i.e.} it holds for $(D_1,D_2,D_3)$ as well as for its dual configuration.

\begin{thm} \label{atmostfour}
There exist at most four minimal disks bounded by $(D_1,D_2,D_3)$ or
its dual configuration and with helicoidal ends of parameters
$(A,\alpha)$, $(B,\beta)$ and $(C,\gamma)$ at $0$, $\infty$ and $1$
respectively. 
\end{thm}

\begin{proof}
By lemmas \ref{fourschwarzian} and \ref{fourschwarziandoubleroot}, it
suffices to prove that for each 
possibility of the Schwarzian derivative there exists at most one minimal
immersion bounded by $(D_1,D_2,D_3)$ or its dual configuration and
with helicoidal ends of parameters $(A,\alpha)$, $(B,\beta)$, $(C,\gamma)$.

Assume that the function $\Theta$ is known. Then the set of the
solutions of equation (\ref{eqk}) on $\Sigma$ is a vector space
generated by two independent solutions. Thus, if $(g,\omega)$ and
$(\tilde g,\tilde\omega)$
are the Weierstrass data of two minimal immersions corresponding to
$\Theta$, then $g$ and $\tilde g$ are quotients of linear combinations
of these two independent solutions, so there
exists a M\"obius transform $\mu:\bar\C\to\bar\C$ such that $\tilde
g=\mu\circ g$. But the value of the Gauss map at each end is uniquely
determined, so we have $g(0)=\tilde g(0)$, $g(1)=\tilde g(1)$ and
$g(\infty)=\tilde g(\infty)$. Moreover, $g(0)$, $g(1)$ and $g(\infty)$
are pairwise distinct (since the straight lines do not lie in parallel
planes), so $\mu$ is the identity, and so $\tilde g=g$ and
$\tilde\omega=\omega$. 
\end{proof}

\subsection{Some facts about the hypergeometric differential equation}
\label{hypergeometricequation}

In this section we recall some facts about the hypergeometric
differential equation and hypergeometric series that will be useful to
give explicit examples of minimal disks bounded by three straight lines.

Let $s_1$, $s_2$ and $s_3$ be three complex numbers such that
$s_3\notin-\N^*$. The hypergeometric series is defined by
$$\hg(s_1,s_2;s_3;z)=\frac{\Gamma(s_3)}{\Gamma(s_1)\Gamma(s_2)}
\sum_{n=0}^{\infty}
\frac{\Gamma(s_1+n)\Gamma(s_2+n)}{\Gamma(s_3+n)}z^n$$
for $|z|<1$. 

We use the notations of section \ref{hopfandspinors}.
We define eight numbers
\begin{equation} \label{eight}
s_{\pm\pm\pm}=\frac{1\pm\alpha\pm\beta\pm\gamma}{2}.
\end{equation}
These numbers are noninteger.
We also set
\begin{equation} \label{Pi}
\begin{array}{ccc}
\Pi & = & (1+\alpha+\beta+\gamma)(1-\alpha+\beta+\gamma)
(1+\alpha-\beta+\gamma) \\ 
& & \times(1+\alpha+\beta-\gamma)(1-\alpha-\beta+\gamma)
(1-\alpha+\beta-\gamma) \\
& & \times(1+\alpha-\beta-\gamma)(1-\alpha-\beta-\gamma).
\end{array}
\end{equation}
We consider the following hypergeometric equation on $\Sigma$:
\begin{equation} \label{he}
w''+\left(\frac{1-\alpha}{z}+\frac{1-\gamma}{z-1}\right)w'
+s_{---}s_{-+-}\frac{w}{z(z-1)}=0.
\end{equation}
This equation is usually denoted by 
$\mathrm{P}\left\{
\begin{array}{ccc}
0 & s_{---} & 0 \\
\alpha & s_{-+-} & \gamma
\end{array}
z\right\}$.
Since we have $\alpha,\beta,\gamma\notin\Z$, the fundamental system of linear
independent solutions of hypergeometric equation (\ref{he}) at the
singular points is given (see section 2.2, paragraph 1 in \cite{mos},
or \cite{bobenko})  
\begin{itemize}
\item on $\Sigma_0$ by
$$w_1^{(0)}(z)=\hg(s_{---},s_{-+-};1-\alpha;z),$$
$$w_2^{(0)}(z)=z^{\alpha}\hg(s_{+--},s_{++-};1+\alpha;z),$$
\item on $\Sigma_1$ by
\begin{eqnarray*}
w_1^{(1)}(z) & = & \hg(s_{---},s_{-+-};1-\gamma;1-z) \\
& = & z^{\alpha}\hg(s_{+--},s_{++-};1-\gamma;1-z),
\end{eqnarray*}
\begin{eqnarray*}
w_2^{(1)}(z) & = & (1-z)^{\gamma}\hg(s_{--+},s_{-++};1+\gamma;1-z) \\
& = & z^{\alpha}(1-z)^{\gamma}\hg(s_{+-+},s_{+++};1+\gamma;1-z),
\end{eqnarray*}
\item on $\Sigma_\infty$ by
$$w_1^{(\infty)}(z)=z^{-s_{---}}\hg(s_{---},s_{+--};1-\beta;z^{-1}),$$
$$w_2^{(\infty)}(z)=z^{-s_{-+-}}\hg(s_{-+-},s_{++-};1+\beta;z^{-1}).$$
\end{itemize}
The second expressions for $w_1^{(1)}$ and $w_2^{(1)}$ are obtained
using first formula of section 2.4.1 in \cite{mos}.

For $z\in(-1,0)$ we have $w_1^{(0)}(z)\in\R$ and $w_2^{(0)}(z)\in
e^{i\pi\alpha}\R$; for $z\in(0,1)$ we have $w_1^{(0)}(z)\in\R$ and
$w_2^{(0)}(z)\in\R$.
For $z\in(0,1)$ we have $w_1^{(1)}(z)\in\R$ and $w_2^{(1)}(z)\in\R$;
for $z\in(1,2)$ we have $w_1^{(1)}(z)\in\R$ and
$w_2^{(1)}(z)\in e^{-i\pi\gamma}\R$.
For $z\in(1,+\infty)$ we have $w_1^{(\infty)}(z)\in\R$ and
$w_2^{(\infty)}(z)\in\R$; for $z\in(-\infty,-1)$ we have 
$w_1^{(\infty)}(z)\in e^{-i\pi s_{---}}\R$ and
$w_2^{(\infty)}(z)\in e^{-i\pi s_{-+-}}\R$.

These solutions are connected in the following way. On
$\Sigma_0\cap\Sigma_1$ we have
$$\left(\begin{array}{c}
w_1^{(0)} \\
w_2^{(0)}
\end{array}\right)=\nu\left(\begin{array}{c}
w_1^{(1)} \\
w_2^{(1)}
\end{array}\right)$$ where
\begin{equation} \label{coeffnu}
\nu=\left(\begin{array}{cc}
\nu_{11} & \nu_{12} \\
\nu_{21} & \nu_{22}
\end{array}\right)=\left(\begin{array}{cc}
\frac{\Gamma(1-\alpha)\Gamma(\gamma)}{\Gamma(s_{--+})\Gamma(s_{-++})}
&
\frac{\Gamma(1-\alpha)\Gamma(-\gamma)}{\Gamma(s_{---})\Gamma(s_{-+-})}
\\
\frac{\Gamma(1+\alpha)\Gamma(\gamma)}{\Gamma(s_{+-+})\Gamma(s_{+++})}
&
\frac{\Gamma(1+\alpha)\Gamma(-\gamma)}{\Gamma(s_{+--})\Gamma(s_{++-})}
\end{array}\right)
\end{equation}
(we used fourth formula of section 2.4.1 in \cite{mos} to compute this
matrix).
In the same way we have
$$\left(\begin{array}{c}
w_1^{(0)} \\
w_2^{(0)}
\end{array}\right)=\hat\nu
\left(\begin{array}{c}
w_1^{(\infty)} \\
w_2^{(\infty)}
\end{array}\right)$$ where
\begin{equation} \label{coeffnuhat}
\begin{array}{ccc}
\hat\nu & = & \left(\begin{array}{cc}
\hat\nu_{11} & \hat\nu_{12} \\
\hat\nu_{21} & \hat\nu_{22}
\end{array}\right) \\
& = & \left(\begin{array}{cc}
e^{i\pi s_{---}}
\frac{\Gamma(1-\alpha)\Gamma(\beta)}{\Gamma(s_{-+-})\Gamma(s_{-++})}
& 
e^{i\pi s_{-+-}}
\frac{\Gamma(1-\alpha)\Gamma(-\beta)}{\Gamma(s_{---})\Gamma(s_{--+})}
\\
e^{i\pi s_{+--}}
\frac{\Gamma(1+\alpha)\Gamma(\beta)}{\Gamma(s_{++-})\Gamma(s_{+++})}
&
e^{i\pi s_{++-}}
\frac{\Gamma(1+\alpha)\Gamma(-\beta)}{\Gamma(s_{+--})\Gamma(s_{+-+})}
\end{array}\right)
\end{array}
\end{equation}
(we used fifth formula of section 2.4.1 in \cite{mos} to compute this
matrix; we should notice that in this formula $(-z)^{-a}$ is
defined with $\arg(-z)\in(-\pi,\pi)$, and thus with this convention we
get $(-z)^{-a}=e^{i\pi a}z^{-a}$). This last formula is actually valid
for the analytic continuations of the solutions on $\Sigma_0$ and
$\Sigma_\infty$, since the intersection of these two domains is empty.

In the sequel, $\sigma_1$ and $\sigma_2$ will denote the solutions on
$\Sigma$ of hypergeometric equation (\ref{he}) such that 
$$\sigma_1=w_1^{(0)},\quad\sigma_2=w_2^{(0)}$$ on $\Sigma_0$.

\begin{lemma} \label{lemmasigma1sigma2}
We have $\sigma_1\sigma_2'-\sigma_1'\sigma_2=
\alpha z^{\alpha-1}(1-z)^{\gamma-1}$.
\end{lemma}

\begin{proof}
Since $\sigma_1$ and $\sigma_2$ are solutions of (\ref{he}), we get
that $\sigma_1\sigma_2'-\sigma_1'\sigma_2$ is a solution on $\Sigma$
of the following equation:
$$w'=-\left(\frac{1-\alpha}{z}+\frac{1-\gamma}{z-1}\right)w.$$
Thus it is proportional to $z^{\alpha-1}(1-z)^{\gamma-1}$. And since
$\sigma_1\sigma_2'-\sigma_1'\sigma_2\sim\alpha z^{\alpha-1}$ when
$z\to 0$, we get the announced expression.
\end{proof}

\begin{lemma} \label{lemmacoeffnu}
We have $\nu_{12}\nu_{21}=-t^2\nu_{11}\nu_{22}$ where $t$ has been defined
in section \ref{geometricconfiguration}.
\end{lemma}

\begin{proof}
Using that $\Gamma(1+z)=z\Gamma(z)$ and
$\Gamma(1-z)\Gamma(z)=\frac{\pi}{\sin(\pi z)}$ (see section 1.1 in
\cite{mos}), we compute that
\begin{eqnarray*}
\nu_{12}\nu_{21} & = & \frac
{\Gamma(1-\alpha)\Gamma(-\gamma)\Gamma(1+\alpha)\Gamma(\gamma)}
{\Gamma(s_{---})\Gamma(s_{-+-})\Gamma(s_{+-+})\Gamma(s_{+++})} \\
& = & -\frac{\alpha}{\gamma}
\frac{\Gamma(1-\alpha)\Gamma(1-\gamma)\Gamma(\alpha)\Gamma(\gamma)}
{\Gamma(s_{---})\Gamma(s_{-+-})\Gamma(1-s_{-+-})\Gamma(1-s_{---})} \\
& = & -\frac{\alpha}{\gamma}
\frac{\sin(\pi s_{---})\sin(\pi s_{-+-})}
{\sin(\pi\alpha)\sin(\pi\gamma)},
\end{eqnarray*}
\begin{eqnarray*}
\nu_{11}\nu_{22} & = & \frac
{\Gamma(1-\alpha)\Gamma(\gamma)\Gamma(1+\alpha)\Gamma(-\gamma)}
{\Gamma(s_{--+})\Gamma(s_{-++})\Gamma(s_{+--})\Gamma(s_{++-})} \\
& = & -\frac{\alpha}{\gamma}
\frac{\Gamma(1-\alpha)\Gamma(\gamma)\Gamma(\alpha)\Gamma(1-\gamma)}
{\Gamma(s_{--+})\Gamma(s_{-++})\Gamma(1-s_{-++})\Gamma(1-s_{--+})} \\
& = & -\frac{\alpha}{\gamma}
\frac{\sin(\pi s_{--+})\sin(\pi s_{-++})}
{\sin(\pi\alpha)\sin(\pi\gamma)}.
\end{eqnarray*}

Thus proving that $\nu_{12}\nu_{21}=-t^2\nu_{11}\nu_{22}$ is
equivalent to prove  that
\begin{equation} \label{eqs}
\sin(\pi s_{---})\sin(\pi s_{-+-})=
-t^2\sin(\pi s_{--+})\sin(\pi s_{-++}).
\end{equation}

But we have
\begin{eqnarray*}
2\sin(\pi s_{---})\sin(\pi s_{-+-})=
\cos(\pi\beta)-\cos(\pi(-\alpha-\gamma+1)) \\
\quad=\cos(\pi\beta)+\cos(\pi\alpha)\cos(\pi\gamma)
-\sin(\pi\alpha)\sin(\pi\gamma),
\end{eqnarray*}
\begin{eqnarray*}
2\sin(\pi s_{--+})\sin(\pi s_{-++})=
\cos(\pi\beta)-\cos(\pi(-\alpha+\gamma+1)) \\
\quad=\cos(\pi\beta)+\cos(\pi\alpha)\cos(\pi\gamma)
+\sin(\pi\alpha)\sin(\pi\gamma).
\end{eqnarray*}
Thus, since $\frac{1-t^2}{1+t^2}=\cos\theta$, condition
($\ref{eqs}$) is equivalent to
\begin{equation*}
\cos\theta=\frac{\cos(\pi\beta)+\cos(\pi\alpha)\cos(\pi\gamma)}
{\sin(\pi\alpha)\sin(\pi\gamma)},
\end{equation*}
which is satified according to lemma \ref{costheta} and since
$\pi\alpha$, $\pi\beta$ and $\pi\gamma$ are congruent to
$\pi\alpha_0$, $\pi\beta_0$ and $\pi\gamma_0$ modulo $2\pi$.
\end{proof}

\begin{lemma} \label{lemmacoeffnuhat}
We have $\hat\nu_{12}\hat\nu_{21}=-\hat{t}^2\hat\nu_{11}\hat\nu_{22}$
where $\hat{t}$ has been defined in section \ref{geometricconfiguration}.
\end{lemma}

\begin{proof}
Since $e^{i\pi s_{---}}e^{i\pi s_{++-}}=e^{i\pi(1-\gamma)}=
e^{i\pi s_{-+-}}e^{i\pi s_{+--}}$,
proceeding as in lemma \ref{lemmacoeffnu} and exchanging the roles of
$\beta$ and $\gamma$, we obtain that the equality of the lemma is equivalent to
\begin{equation*}
\cos\hat\theta=\frac{\cos(\pi\gamma)+\cos(\pi\alpha)\cos(\pi\beta)}
{\sin(\pi\alpha)\sin(\pi\beta)}.
\end{equation*}
This condition is satified according to lemma
\ref{costheta} and since
$\pi\alpha$, $\pi\beta$ and $\pi\gamma$ are congruent to
$\pi\alpha_0$, $\pi\beta_0$ and $\pi\gamma_0$ modulo $2\pi$.
\end{proof}

\begin{lemma} \label{lemmanunuhat}
We have $$\frac{\hat{t}\hat\nu_{11}}{\hat\nu_{21}}=
e^{-i\pi\alpha}\frac{t\nu_{11}}{\nu_{21}}.$$
\end{lemma}

\begin{proof}
Proving this equality is equivalent to prove that
$$t\sin(\pi s_{--+})=\hat{t}\sin(\pi s_{-+-}).$$

We have computed that
$$t^2=-\frac{\sin(\pi s_{---})\sin(\pi s_{-+-})}
{\sin(\pi s_{--+})\sin(\pi s_{-++})},\quad
\hat{t}^2=-\frac{\sin(\pi s_{---})\sin(\pi s_{--+})}
{\sin(\pi s_{-+-})\sin(\pi s_{-++})}.$$
We deduce that $t^2\sin^2(\pi s_{--+})=\hat{t}^2\sin^2(\pi
s_{-+-})$. Since $\sin\theta=\frac{2t}{1+t^2}$ and
$\sin\hat\theta=\frac{2\hat t}{1+\hat t^2}$, we deduce from lemma
\ref{signs} that $t$ and $\hat t$ have the same sign. So it now
suffices to prove that $\sin(\pi s_{--+})$ and $\sin(\pi s_{-+-})$
have the same sign.

We have $2\sin(\pi s_{--+})\sin(\pi s_{-+-})
=\cos(\pi(\gamma-\beta))-\cos(\pi(1-\alpha))
=\cos(\pi(\gamma_0-\beta_0))-\cos(\pi(1-\alpha_0))$.
We also have $1-\alpha_0\in(0,1)$ and $\gamma_0-\beta_0\in(-1,1)$, and by
(\ref{sphericalangles}) we have $1-\alpha_0>|\gamma_0-\beta_0|$, so we get
$\cos(\pi(1-\alpha_0))<\cos(\pi(\gamma_0-\beta_0))$ and
$2\sin(\pi s_{--+})\sin(\pi s_{-+-})>0$. This proves the lemma.
\end{proof}

\subsection{Existence of a minimal surface bounded by three lines}

In this section we give explicit examples of minimal disks bounded
by three lines in generic position. We use the notations of section
\ref{hopfandspinors}.

Let $a$, $b$ and $c$ be three real
numbers. For $j=1,2$ and $z\in\Sigma$, we set
\begin{equation} \label{defK}
K_j=z^{\frac{1-\alpha}{2}}(1-z)^{\frac{1-\gamma}{2}}
((a+bz)\sigma_j+cz(1-z)\sigma_j').
\end{equation}
These functions were introduced by Riemann in his memoir
\cite{riemann} (where they are denoted by $k_1$ and $k_2$).

\begin{lemma}[Riemann, \cite{riemann}] \label{lemmariemann}
The functions $K_1$ and $K_2$ satisfy
$$K_1K_2'-K_1'K_2=z^{1-\alpha}(1-z)^{1-\gamma}
(\sigma_1\sigma_2'-\sigma_1'\sigma_2)F(z)$$
with
\begin{equation} \label{eqF}
\begin{array}{ccc}
F(z) & = & a(a+c\alpha)(1-z)+(a+b)(a+b-c\gamma)z \\
& & -(b+s_{---}c)(b+s_{-+-}c)z(1-z).
\end{array}
\end{equation}
\end{lemma}

\begin{proof}
Using the fact that $\sigma_j$ is a solution of (\ref{he}), we compute
that
\begin{eqnarray*}
K_j' & = & z^{-\frac{1+\alpha}{2}}(1-z)^{-\frac{1+\gamma}{2}}\times \\
& & ((n_1z^2+n_2z(1-z)+n_3(1-z)^2)\sigma_j) \\
& & +(n_4z^2(1-z)+n_5z(1-z)^2)\sigma_j')
\end{eqnarray*}
with
$$n_1=-\frac{1-\gamma}{2}(a+b),\quad
n_2=\frac{\gamma-\alpha}{2}a+\frac{3-\alpha}{2}b+cs_{---}s_{-+-},$$
$$n_3=\frac{1-\alpha}{2}a,\quad n_4=a+b-\frac{1+\gamma}{2}c,\quad
n_5=a+\frac{1+\alpha}{2}c.$$
Thus we get
\begin{eqnarray*}
K_1K_2'-K_1'K_2 & = & z^{1-\alpha}(1-z)^{1-\gamma}
(\sigma_1\sigma_2'-\sigma_1'\sigma_2)\times \\
& & (((a+b)n_4-cn_1)z^2 \\
& & +(an_4+(a+b)n_5-cn_2)z(1-z) \\
& & +(an_5-cn_3)(1-z)^2).
\end{eqnarray*}
The last factor in this expression is a polynomial of degree $2$. We
compute that its values at $0$ and $1$ and its degree $2$ coefficient
are the same as those of the polynomial $F$ in the lemma.
\end{proof}

\begin{cor}
The functions $K_1$ and $K_2$ satisfy
\begin{equation} \label{eqK1K2}
K_1K_2'-K_1'K_2=\alpha F, \quad K_1K_2''-K_1''K_2=\alpha F',
\end{equation}
with $F$ as in (\ref{eqF}).
\end{cor}

\begin{proof}
The first formula comes from lemmas \ref{lemmasigma1sigma2} and
\ref{lemmariemann}. We obtain the second one by differentiation.
\end{proof}

\begin{lemma} \label{mureal}
Let $\lambda_1,\lambda_2,\mu_1,\mu_2$ be complex numbers such that
$\lambda_1^{-1}\mu_2-\lambda_2^{-1}\mu_1\neq 0$. 
Then the spinors $k_1=\lambda_1^{-1}K_1+\lambda_2^{-1}K_2$ and
$k_2=\mu_1K_1+\mu_2K_2$ define a conformal minimal immersion
$x:\Sigma\to\R^3$, with possibly a singular point at the root of $F$
when $F$ has a double root.
The immersion $x$ has, up to a translation in $\R^3$, a helicoidal end
bounded by $D_1$ and $D_2$ and of parameters $(A,\alpha)$ if and only if
$\lambda_2=\mu_1=0$, $\lambda_1\in i\R^*$, $\mu_2\in\R^*$ and 
$\alpha\lambda_1^{-1}\mu_2a(a+\alpha c)=i\frac{A\alpha}{2\pi}$.
\end{lemma}

\begin{proof}
The singularities of $x$ correspond to the common zeros of $k_1$ and
$k_2$, and thus to the common zeros of $K_1$ and $K_2$. Using
(\ref{eqK1K2}) we conclude that a singular point of $x$ is necessarily
a double root of $F$.

Let $(g,\omega)$ be the Weierstrass data of $x$. We will use the
expressions of $\sigma_1$ and $\sigma_2$ 
valid in $\Sigma_0$, that is $\sigma_1=w_1^{(0)}$ and
$\sigma_2=w_2^{(0)}$.  We have
$$g=\frac{\mu_1K_1+\mu_2K_2}{\lambda_1^{-1}K_1+\lambda_2^{-1}K_2},$$
$$k_1k_2'-k_1'k_2=
(\lambda_1^{-1}\mu_2-\lambda_2^{-1}\mu_1)(K_1K_2'-K_1'K_2)
=\alpha(\lambda_1^{-1}\mu_2-\lambda_2^{-1}\mu_1)F$$ with $F$ as in
(\ref{eqF}).

Assume that $x$ has, up to a translation in $\R^3$, a helicoidal end
bounded by $D_1$ and $D_2$ and of parameters $(A,\alpha)$. Then when
$z\to 0$ we have $Q\sim i\frac{A\alpha}{2\pi}z^{-2}\rmd z^2$, so we get
$i\frac{A\alpha}{2\pi}=\alpha(\lambda_1^{-1}\mu_2-\lambda_2^{-1}\mu_1)F(0)
=\alpha(\lambda_1^{-1}\mu_2-\lambda_2^{-1}\mu_1)a(a+\alpha c)$. This
implies in particular that 
$\lambda_1^{-1}\mu_2-\lambda_2^{-1}\mu_1\in i\R$,
$a\neq0$ and $a+\alpha c\neq0$.

Then we have $K_1(z)\sim az^{\frac{1-\alpha}{2}}$ and
$K_2(z)\sim(a+\alpha c)z^{\frac{1+\alpha}{2}}$ when $z\to 0$.
We must have $g(z)\sim i\rho z^\alpha$ for some $\rho\in\R^*$ (see the
proof of lemma \ref{parameters}), since $\alpha\in(0,1)+2\Z$. This
implies that $\mu_1=0$ if $\alpha>0$ and $\lambda_2=0$ if $\alpha<0$.

We deal with the case where $\alpha>0$. In this case we have 
$g=\frac{\mu_2K_2}{\lambda_1^{-1}K_1+\lambda_2^{-1}K_2}
\sim\lambda_1\mu_2\frac{a+\alpha c}{a}z^\alpha$,
so $\lambda_1\mu_2\in i\R$. And since we also have
$\lambda_1^{-1}\mu_2\in i\R$ (because $\lambda_2^{-1}\mu_1=0$), we get
$\mu_2^2\in\R$ and $\lambda_1^2\in\R$. Moreover, we have $g(z)\in
i\bar\R$ if $z\in(0,1)$ (since $D_2$ is the $x_1$-axis), and the
functions $K_1$ and $K_2$ take real values on $(0,1)$, so we also have
$\lambda_2\mu_2\in i\R$, and so $\lambda_2^2\in\R$.
Since $x$ maps
$(-1,0)$ onto a straight line, and since
$\rmd(x_1+ix_2)=\bar\omega-g^2\omega
=\overline{z^{-2}(z-1)^{-2}k_1^2\rmd z}
-z^{-2}(z-1)^{-2}k_2^2\rmd z$, the argument of $\overline{k_1^2}-k_2^2$
must be constant on $(-1,0)$.
We have $\overline{k_1^2}-k_2^2
%\overline{(\lambda_1^{-1}K_1+\lambda_2^{-1}K_2)^2}-\mu_2^2K_2^2$.
=\lambda_1^{-2}\overline{K_1^2}
+2\overline{\lambda_1^{-1}\lambda_2^{-1}K_1K_2}
+\lambda_2^{-2}\overline{K_2^2}-\mu_2^2K_2^2$; on the other hand we have
$\arg{K_1(z)}\equiv\pi\frac{1-\alpha}{2}\mod\pi$ and
$\arg{K_2(z)}\equiv\pi\frac{1+\alpha}{2}\mod\pi$ 
if $z\in(-1,0)$, so we conclude that $\lambda_2=0$.
Thus we have $\omega\sim\lambda_1^{-2}a^2z^{-1-\alpha}\rmd z$ and
$g^2\omega\sim\mu_2^2(a+\alpha c)^2z^{-1+\alpha}\rmd z$.
Since $x_1\to+\infty$ when $z\to 0$ with $z$ real and positive (by
condition 2 in definition \ref{helicoidalend}), we get $\lambda_1^2<0$
and $\mu_2^2>0$, which implies $\lambda_1\in i\R^*$ and
$\mu_2\in\R^*$.

We proceed in the same way in the case where $\alpha<0$.

Conversely, assume that $\lambda_2=\mu_1=0$, $\lambda_1\in i\R^*$,
$\mu_2\in\R^*$ and
$\alpha\lambda_1^{-1}\mu_2a(a+\alpha c)=i\frac{A\alpha}{2\pi}$.
Then we have $\arg{\overline{k_1^2}}=\arg(-k_2^2)=\pi\alpha$ on
$(-1,0)$ and $\arg{\overline{k_1^2}}=\arg(-k_2^2)=\pi$ on $(0,1)$, and
we have $\rmd(x_1+ix_2)=\bar\omega-g^2\omega\sim
\lambda_1^{-2}a^2\overline{z^{-1-\alpha}\rmd z}
-\mu_2^2(a+\alpha c)^2z^{-1+\alpha}\rmd z$. This proves that $x$ has
an end bounded by two lines $D_1'$ and $D_2'$ that are parallel to
$D_1$ and $D_2$ respectively. Moreover, we have
$g(z)\sim\lambda_1\mu_2\frac{a+\alpha c}{a}z^\alpha$ 
and $Q\sim\alpha\lambda_1^{-1}\mu_2a(a+\alpha c)z^{-2}\rmd z^2
=i\frac{A\alpha}{2\pi}z^{-2}\rmd z^2$, so $x$ has a helicoidal end of
parameters $(A,\alpha)$ by lemma \ref{helicoidalend}.
This finally implies that $\rmD(D_1',D_2')=-A=\rmD(D_1,D_2)$, and so
$D_1'$ and $D_2'$ are the images of $D_1$ and $D_2$ by a translation
in $\R^3$.
\end{proof}

In the sequel we will study the minimal immersion $x:\Sigma\to\R^3$
(with possibly a singularity at a double root of $F$)
given by the spinors 
$k_1=\lambda^{-1}K_1$ and $k_2=\mu K_2$ for some
$\lambda\in\C^*$ and some $\mu\in\C^*$. We first notice the following fact.

\begin{rem} \label{remarklambdamu}
For $\rho\in\R^*$, the transformation
$$(a,b,c,\lambda,\mu)\to (\rho a,\rho b,\rho
c,\rho\lambda,\rho^{-1}\mu)$$
does not change the 
Weierstrass data $(g,\omega)$ (and consequently does not change the
immersion $x$), and changes $\lambda^{-1}\mu$ into
$\rho^{-2}\lambda^{-1}\mu$. Thus, without loss of generality, we can
assume that $|\lambda\mu^{-1}|=|\alpha|$.
\end{rem}

\begin{rem}
Replacing $(\lambda,\mu)$ by $(-i\lambda,i\mu)$ would change $g$ into
$-g$ and $\omega$ into $-\omega$. Thus the immersion $x$ would be
replaced by its image by the reflexion about the $x_3$-axis.
\end{rem}

\begin{prop} \label{threelines}
Let $\lambda$ and $\mu$ be two nonzero complex numbers such that
$|\lambda\mu^{-1}|=|\alpha|$.
Let $x:\Sigma\to\R^3$ be the conformal minimal immersion (with
possibly a singularity at a double root of $F$) whose
Weierstrass data are given by
$$g=\frac{k_2}{k_1},\quad\omega=z^{-2}(z-1)^{-2}k_1^2,$$
where $$k_1=\lambda^{-1}K_1,\quad k_2=\mu K_2$$ with $K_1$ and $K_2$
as in formula (\ref{defK}).

Then, up to a translation, the immersion $x$ maps $(-\infty,0)$,
$(0,1)$ and $(1,+\infty)$ to $D_1$, $D_2$ 
and $D_3$ respectively, and has helicoidal ends of parameters
$(A,\alpha)$, $(B,\beta)$ and $(C,\gamma)$ at $0$, $\infty$ and $1$
respectively if and only if the following conditions hold:
\begin{enumerate}
\item[1.] $\lambda=-\varepsilon i\alpha\mu$ where 
$\varepsilon$ is the sign of
$\varepsilon_0\frac{\nu_{11}}{\alpha\nu_{21}}$ with $\nu$ as defined
in section \ref{hypergeometricequation},
\item[2.] $\varepsilon\alpha\mu^2=\frac{t\nu_{11}}{\nu_{21}}$ (this
implies in particular that $\mu$ is real),
\item[3.] $F=\varepsilon\varphi$ with $\varphi$ as in (\ref{eqphi}) and
$F$ as in (\ref{eqF}), {\it i.e.}
the real numbers $a$, $b$ and $c$ satisfy
\begin{eqnarray} \label{systemabc}
\left\{\begin{array}{ccc}
\varepsilon\frac{A\alpha}{2\pi} & = & a(a+\alpha c) \\
\varepsilon\frac{B\beta}{2\pi} & = & (b+s_{---}c)(b+s_{-+-}c) \\
\varepsilon\frac{C\gamma}{2\pi} & = & (a+b)(a+b-\gamma c).
\end{array}\right.
\end{eqnarray}
\end{enumerate}
\end{prop}

\begin{proof}
We denote by $\Delta_1$,
$\Delta_2$, $\Delta_3$, $\Delta_1^{\infty}$, $\Delta_1^0$,
$\Delta_3^0$ and $\Delta_3^{\infty}$ the images by $x$ of
$(-\infty,0)$, $(0,1)$, $(1,\infty)$, $(-\infty,-1)$, $(-1,0)$,
$(1,2)$ and $(2,+\infty)$.

If $x$ has, up to a translation in $\R^3$, a helicoidal end at $0$
bounded by $D_1$ and $D_2$ and of parameters $(A,\alpha)$, then by
lemma \ref{mureal} we have $\lambda\in i\R^*$, $\mu\in\R^*$ and
$\alpha\lambda^{-1}\mu a(a+\alpha c)=i\frac{A\alpha}{2\pi}$. Since
$|\lambda\mu^{-1}|=|\alpha|$, we have $\lambda=-\varepsilon i\alpha\mu$
with $\varepsilon=\pm 1$, and so the first equality in
(\ref{systemabc}) holds. 

From now on we assume that
\begin{equation} \label{choicelambdamu}
\lambda=-\varepsilon i\alpha\mu
\end{equation}
with $\varepsilon=\pm 1$ and $\mu\in\R^*$, and that
$\frac{A\alpha}{2\pi}=a(a+\alpha c)$. Then by lemma \ref{mureal} we
can assume that the immersion $x$ has a helicoidal end at $0$ 
bounded by $D_1$ and $D_2$ and of parameters $(A,\alpha)$ (it suffices
to consider the good translation in $\R^3$). Moreover we have
$\Delta_1^0\subset D_1$, $\Delta_2\subset D_2$, $\Delta_1^0$ contains
a half of $D_1$ in the direction of $\vec{u}_1$, and $\Delta_2$ contains
a half of $D_2$ in the direction of $-\vec{u}_2$.

We now prove the necessity of (\ref{systemabc}). By (\ref{eqK1K2})
we have $k_1'k_2-k_1'k_2=\lambda^{-1}\mu\alpha F=\varepsilon iF$ with
$F$ as in (\ref{eqF}). On the other hand, if $x$ has helicoidal ends
of parameters $(A,\alpha)$, $(B,\beta)$ and $(C,\gamma)$ at $0$,
$\infty$ and $1$ respectively, then equation (\ref{eqk1k2}) holds with
$\varphi$ as in (\ref{eqphi}), so we get $\varphi=\varepsilon F$,
which implies that (\ref{systemabc}) holds.

From now on we assume that $a$, $b$, $c$ and $\varepsilon$ satisfy
(\ref{systemabc}). We study the behaviour of $x$ at $z=1$.

The rotation $R$ defined in section \ref{geometricconfiguration} moves
$D_3$ onto a horizontal line oriented by the vector
$(\cos(\pi\gamma),-\sin(\pi\gamma),0)$, the vector $\vec{v}_1$
to the vector $-\vec{e}_3$, and it does not change $D_2$. Let
$\tilde{x}=R\circ x$. Let 
$(\tilde{g},\tilde{\omega})$ be its Weierstrass data, and
$\tilde{N}$ its Gauss map. There exists a matrix
$$h=\left(\begin{array}{cc}
h_{11} & h_{12} \\
h_{21} & h_{22}
\end{array}\right)\in\mathrm{SU}_2(\C)$$
such that
$$\tilde{g}=\frac{h_{22}g+h_{21}}{h_{12}g+h_{11}},\quad
\tilde{\omega}=(h_{12}g+h_{11})^2\omega.$$
Then the associated spinors can be choosen as
$$\tilde{k}_1=(h_{12}g+h_{11})k_1=(h_{12}k_2+h_{11}k_1),\quad
\tilde{k}_2=\tilde{g}\tilde{k}_1=(h_{22}k_2+h_{21}k_1).$$
We compute that $$h=\frac{1}{\sqrt{1+t^2}}
\left(\begin{array}{cc}
1 & it \\
it & 1
\end{array}\right).$$
Consequently we have
$$\tilde{k}_1=z^{\frac{1-\alpha}{2}}(1-z)^{\frac{1-\gamma}{2}}
((a+bz)\tilde{\sigma}_1+cz(1-z)\tilde{\sigma}_1'),$$
$$\tilde{k}_2=z^{\frac{1-\alpha}{2}}(1-z)^{\frac{1-\gamma}{2}}
((a+bz)\tilde{\sigma}_2+cz(1-z)\tilde{\sigma}_2'),$$
where
$$\left(\begin{array}{c}
\tilde{\sigma}_1 \\
\tilde{\sigma}_2
\end{array}\right)=m
\left(\begin{array}{c}
w_1^{(1)} \\
w_2^{(1)}
\end{array}\right)$$
with
$$m=\left(\begin{array}{cc}
m_{11} & m_{12} \\
m_{21} & m_{22}
\end{array}\right)=\frac{1}{\sqrt{1+t^2}}
\left(\begin{array}{cc}
\lambda^{-1} & \mu it \\
\lambda^{-1}it & \mu
\end{array}\right)
\left(\begin{array}{cc}
\nu_{11} & \nu_{12} \\
\nu_{21} & \nu_{22}
\end{array}\right).$$
These expressions are valid for $z\in\Sigma_1$. We notice that
$m_{11}\in i\R$ and $m_{22}\in\R$ since $\lambda\in i\R$ and
$\mu\in\R$.

We claim that $\tilde{x}$ has a helicoidal end at $z=1$ of parameters
$(C,\gamma)$ bounded by $D_2=R(D_2)$ and $R(D_3)$ if and only if
$m_{12}=m_{21}=0$. The proof of this claim is similar to that of lemma
\ref{mureal}, so we will only outline the proof. We already know that
$\tilde x$ maps $(0,1)$ to a part of $D_2$.

Assume that $\tilde{x}$ has a helicoidal end of parameters
$(C,\gamma)$  bounded by $D_2$ and $R(D_3)$. Then we must have
$g(z)\sim i\rho(1-z)^\gamma$ for some $\rho\in\R^*$. This implies
$m_{21}=0$ if $\gamma>0$ and $m_{12}=0$ if $\gamma<0$ (this follows
from the fact that $w_1^{(1)}(1)=1$ and $w_2^{(1)}(z)\sim(1-z)^\gamma$
when $z\to 1$). Moreover, $\tilde{x}$ maps $(1,2)$ onto a straight
line, so the argument of $\overline{\tilde{k}_1^2}-\tilde{k}_2^2$ is
constant on $(1,2)$. This implies that $m_{12}=0$ if $m_{21}=0$ and
$m_{21}=0$ if  $m_{12}=0$ (because $w_1^{(1)}$, $w_2^{(1)}$ and their
derivatives take real values on $(1,2)$). 

Conversely, we assume that $m_{12}=m_{21}=0$. Since $m_{11}\in i\R$
and $m_{22}\in\R$ we have $\arg{\overline{\tilde k_1^2}}
=\arg(-\tilde k_2^2)=\pi$ on $(1,0)$ and $\arg{\overline{\tilde k_1^2}}
=\arg(-\tilde k_2^2)=-\pi\gamma$ on $(1,2)$, and we have
$\rmd(\tilde x_1+i\tilde x_2)=\overline{\tilde\omega}-\tilde g^2\tilde\omega
\sim m_{11}^2(a+b)^2\overline{(1-z)^{-1-\gamma}\rmd z}
-m_{22}^2(a+b-\gamma c)^2(1-z)^{-1+\gamma}\rmd z$. This proves that
$\tilde x$ has an end bounded by two lines that are parallel to $D_2$
and $R(D_3)$, and the signed distance of these lines is equal to $-C$
because of the third equation in (\ref{systemabc}). Thus, since we
know that $\tilde x$ maps $(0,1)$ to a part of $D_2$, we have proved
that $\tilde{x}$ has a helicoidal end of parameters
$(C,\gamma)$ bounded by $D_2$ and $R(D_3)$. This completes proving the
claim.

The condition $m_{12}=m_{21}=0$ is satisfied if and only
if $\lambda^{-1}\nu_{12}+\mu it\nu_{22}
=\lambda^{-1}it\nu_{11}+\mu\nu_{21}=0$, that is, because of
(\ref{choicelambdamu}), if and only if 
\begin{equation} \label{eqlambda}
\varepsilon\alpha\mu^2=-\frac{\nu_{12}}{t\nu_{22}}=\frac{t\nu_{11}}{\nu_{21}}.
\end{equation}
We recall that $\mu$ must be real.
Thus there exist $\mu\in\R^*$ and $\varepsilon\in\{1,-1\}$ satisfying
(\ref{eqlambda}) if and only if 
$\nu_{12}\nu_{21}=-t^2\nu_{11}\nu_{22}$, which is true by lemma
\ref{lemmacoeffnu}.

Henceforth we assume that $\mu$ and $\varepsilon$ are given by
(\ref{eqlambda}). The real number $\mu$ is defined uniquely up to its
sign, and $\varepsilon$ is the sign of
$\frac{t\nu_{11}}{\alpha\nu_{21}}$, {\it i.e.} of
$\varepsilon_0\frac{\nu_{11}}{\alpha\nu_{21}}$ by lemma
\ref{signs}. Thus the number $\lambda$ is also 
defined uniquely up to its sign by (\ref{choicelambdamu}). Thus we
have proved the necessity of conditions 1 and 2. 

We have $\Delta_1^0\subset D_1$, $\Delta_2=D_2$, $\Delta_3^1\subset
D_3$, $\Delta_1^0$ contains a half of $D_1$ in the direction of
$\vec{u}_1$ and $\Delta_3^1$ contains a half of $D_3$ in the direction
of $-\vec{u}_3$. It now suffices to check that $x$ has a helicoidal
end at $\infty$ bounded by $D_3$ and $D_1$ and of parameters
$(B,\beta)$. 

The isometry $T\circ\hat R$ defined in section \ref{geometricconfiguration}
moves $D_1$ onto the $x_1$-axis oriented by $\vec{e}_1$, $D_3$ onto a
horizontal line oriented by the vector
$(\cos(\pi\beta),-\sin(\pi\beta),0)$ and 
the vector $\vec{v}_{\infty}$ onto the vector 
$-\vec{e}_3$. Let $\hat x=T\circ\hat R\circ x$. Let $(\hat
g,\hat\omega)$ be its Weierstrass data, and $\hat N$ its Gauss map.
There exists a matrix
$$\hat{h}=\left(\begin{array}{cc}
\hat{h}_{11} & \hat{h}_{12} \\
\hat{h}_{21} & \hat{h}_{22}
\end{array}\right)\in\mathrm{SU}_2(\C)$$
such that
$$\hat{g}=\frac{\hat{h}_{22}g+\hat{h}_{21}}{\hat{h}_{12}g+\hat{h}_{11}},\quad
\hat{\omega}=(\hat{h}_{12}g+\hat{h}_{11})^2\omega.$$
Then the associated spinors can be choosen as
$$\hat{k}_1=(\hat{h}_{12}g+\hat{h}_{11})k_1
=(\hat{h}_{12}k_2+\hat{h}_{11}k_1),\quad
\hat{k}_2=\hat{g}\hat{k}_1=(\hat{h}_{22}k_2+\hat{h}_{21}k_1).$$
We compute that $$\hat{h}=\frac{1}{\sqrt{1+\hat{t}^2}}
\left(\begin{array}{cc}
e^{i\pi\frac{\alpha}{2}} & ie^{-i\pi\frac{\alpha}{2}}\hat{t} \\
ie^{i\pi\frac{\alpha}{2}}\hat{t} & e^{-i\pi\frac{\alpha}{2}}
\end{array}\right).$$
Consequently we have 
$$\hat{k}_1=z^{\frac{1-\alpha}{2}}(1-z)^{\frac{1-\gamma}{2}}
((a+bz)\hat{\sigma}_1+cz(1-z)\hat{\sigma}_1'),$$
$$\hat{k}_2=z^{\frac{1-\alpha}{2}}(1-z)^{\frac{1-\gamma}{2}}
((a+bz)\hat{\sigma}_2+cz(1-z)\hat{\sigma}_2'),$$
where
$$\left(\begin{array}{c}
\hat{\sigma}_1 \\
\hat{\sigma}_2
\end{array}\right)=\left(\begin{array}{cc}
\hat{m}_{11} & \hat{m}_{12} \\
\hat{m}_{21} & \hat{m}_{22}
\end{array}\right)
\left(\begin{array}{c}
w_1^{(\infty)} \\
w_2^{(\infty)}
\end{array}\right)$$
with
$$\left(\begin{array}{cc}
\hat{m}_{11} & \hat{m}_{12} \\
\hat{m}_{21} & \hat{m}_{22}
\end{array}\right)=\frac{1}{\sqrt{1+\hat{t}^2}}
\left(\begin{array}{cc}
\lambda^{-1}e^{i\pi\frac{\alpha}{2}} & 
\mu ie^{-i\pi\frac{\alpha}{2}}\hat{t} \\
\lambda^{-1}ie^{i\pi\frac{\alpha}{2}}\hat{t} &
\mu e^{-i\pi\frac{\alpha}{2}}
\end{array}\right)
\left(\begin{array}{cc}
\hat\nu_{11} & \hat\nu_{12} \\
\hat\nu_{21} & \hat\nu_{22}
\end{array}\right).$$
These expressions are valid for $z\in\Sigma_\infty$.

We first prove that $\hat m_{12}=\hat m_{21}=0$. Indeed, this
condition is equivalent to 
$\lambda^{-1}e^{i\pi\frac{\alpha}{2}}\hat\nu_{12}
+\mu ie^{-i\pi\frac{\alpha}{2}}\hat{t}\hat\nu_{22}
=\lambda^{-1}ie^{i\pi\frac{\alpha}{2}}\hat{t}\hat\nu_{11}
+\mu e^{-i\pi\frac{\alpha}{2}}\hat\nu_{21}=0$, that is, because of
(\ref{choicelambdamu}), to
\begin{equation} \label{eqlambdahat}
\varepsilon e^{-i\pi\alpha}\alpha\mu^2=
-\frac{\hat\nu_{12}}{\hat{t}\hat\nu_{22}}
=\frac{\hat{t}\hat\nu_{11}}{\hat\nu_{21}}.
\end{equation}
Since $\mu$ is defined by (\ref{eqlambda}), this condition is
equivalent to
$$\hat\nu_{12}\hat\nu_{21}=-\hat{t}^2\hat\nu_{11}\hat\nu_{22},\quad
\frac{\hat{t}\hat\nu_{11}}{\hat\nu_{21}}=
e^{-i\pi\alpha}\frac{t\nu_{11}}{\nu_{21}},$$
which holds by lemmas \ref{lemmacoeffnuhat} and
\ref{lemmanunuhat}. This completes proving that $\hat m_{12}=\hat m_{21}=0$.

We claim that $\hat x$ has a helicoidal end at $z=\infty$ of parameters
$(B,\beta)$ bounded by $T\circ\hat R(D_3)$ and $T\circ\hat R(D_1)$.
Since $\lambda\in i\R$
and $\mu\in\R$, we have $\arg(\hat m_{11}^2)=-\pi(\beta+\gamma)$ and
$\arg(\hat m_{22}^2)=\pi(1+\beta-\gamma)$, and so
$\arg{\overline{\hat k_1^2}}=\arg(-\hat k_2^2)=\pi(1+\beta)$ on
$(1,+\infty)$ and 
$\arg{\overline{\hat k_1^2}}=\arg(-\hat k_2^2)=0$ on
$(-\infty,-1)$. We also have $\rmd(\hat x_1+i\hat x_2)
=\overline{\hat\omega}-\hat g^2\hat\omega
\sim\kappa_1\overline{z^{-1+\beta}\rmd z}
+\kappa_2z^{-1-\beta}\rmd z$ for some $\kappa_1,\kappa_2\in\C^*$.
This proves that $\hat x$ has an end bounded by two lines that are
parallel to $T\circ\hat R(D_3)$ and $T\circ\hat R(D_3)$, and the
signed distance of these lines is equal to $-B$ 
because of the second equation in (\ref{systemabc}). Thus, since we
know that $\hat x$ maps $(1,2)$ to a part of $T\circ\hat R(D_3)$, we
have proved that $\hat x$ has a helicoidal end of parameters
$(B,\beta)$ bounded by $T\circ\hat R(D_3)$ and $T\circ\hat
R(D_3)$. This completes proving the claim, and this also prove that
$\Delta_1=D_1$ and $\Delta_3=D_3$. 

Thus the immersion $x$ is bounded by $D_1$, $D_2$ and $D_3$, and it
has helicoidal ends of parameters $(A,\alpha)$, $(B,\beta)$ and
$(C,\gamma)$ at $0$, $\infty$ and $1$ respectively.
\end{proof}

In his memoir \cite{riemann}, Riemann proved the necessity of
(\ref{systemabc}), and only with $\varepsilon=1$: this comes from the fact
that he did not take orientations precisely in consideration. Riemann
also noticed the following fact (page 335).

\begin{lemma}[Riemann, \cite{riemann}] \label{equivalence}
System (\ref{systemabc}) is equivalent to
\begin{eqnarray} \label{systempqr}
\left\{\begin{array}{ccc}
p^2-\alpha^2(p+q+r)^2 & = & \varepsilon\frac{A\alpha}{2\pi} \\
q^2-\beta^2(p+q+r)^2 & = & \varepsilon\frac{B\beta}{2\pi} \\
r^2-\gamma^2(p+q+r)^2 & = & \varepsilon\frac{C\gamma}{2\pi}
\end{array}\right.
\end{eqnarray}
with $p=a+\frac{\alpha}{2}c$, $q=b+\frac{1-\alpha-\gamma}{2}c$ and
$r=-a-b+\frac{\gamma}{2}c$.
\end{lemma}

\begin{thm} \label{theoremsystem}
A real solution $(p,q,r)$ of system (\ref{systempqr}), where
$\varepsilon$ is the sign of
$\varepsilon_0\frac{\nu_{11}}{\alpha\nu_{21}}$, gives a minimal
immersion $x:\Sigma\to\R^3$ (possibly with a singular point at a
double root of $\varphi$ defined by (\ref{eqphi})) bounded by
$(D_1,D_2,D_3)$ and having helicoidal ends of parameters $(A,\alpha)$,
$(B,\beta)$ and $(C,\gamma)$ at $0$, $\infty$ and $1$ respectively. We
will denote it $\cI(\alpha,\gamma,\beta,\varepsilon_0,p,r,q)$.

Moreover, two different real solutions of system (\ref{systempqr})
give the same immersion if and only if they are opposite one to the
other.
\end{thm}

\begin{proof}
The first assertion is a consequence of proposition \ref{mureal} and
lemma \ref{equivalence}.

Assume that $(p,q,r)$ and $(\hat p,\hat q,\hat r)$ are two real
solution of (\ref{systempqr}). Then they correspond to two solutions
$(a,b,c)$ and $(\hat a,\hat b,\hat c)$ of (\ref{systemabc}), which
define functions $K_1$, $K_2$, $\hat K_1$, $\hat K_2$ by
(\ref{defK}). If they define the same immersion, then their
Weierstrass data are equal, and so $\hat K_1=\pm K_1$ and $\hat
K_2=\pm K_2$ (because $\lambda$ and $\mu$ are determined), which
implies $(\hat a,\hat b,\hat c)=\pm(a,b,c)$, and finally $(\hat p,\hat
q,\hat r)=\pm(p,q,r)$. The converse is clear.
\end{proof}

\begin{cor} \label{existencesmallangles}
If $\alpha,\beta,\gamma\in(0,1)$ (that is if $\alpha=\alpha_0$,
$\beta=\beta_0$ and $\gamma=\gamma_0$), then there exists a minimal
immersion $x:\Sigma\to\R^3$ (possibly with a singular point at a 
double root of $\varphi$) bounded by
$(D_1,D_2,D_3)$ and having helicoidal ends of parameters $(A,\alpha)$,
$(B,\beta)$ and $(C,\gamma)$ at $0$, $\infty$ and $1$ respectively.
\end{cor}

\begin{proof}
It suffices to prove that system (\ref{systempqr}) has at least one
real solution. We set
$$p(y)=\sqrt{\varepsilon\frac{A\alpha}{2\pi}+\alpha^2y^2},$$ 
$$q(y)=\sqrt{\varepsilon\frac{B\beta}{2\pi}+\beta^2y^2},$$
$$r(y)=\sqrt{\varepsilon\frac{C\gamma}{2\pi}+\gamma^2y^2},$$
and we define eight functions
$$F_{\pm\pm\pm}(y)=\pm p(y)\pm q(y)\pm r(y)-y$$
for real numbers $y$ such that these expressions are defined. Then
system (\ref{systempqr}) has a real solution if and only if at least one of
the eight equations $F_{\pm\pm\pm}(y)=0$ has a solution $y\in\R$.
When $y\to+\infty$, we have 
$F_{\pm\pm\pm}(y)\sim(\pm\alpha\pm\beta\pm\gamma-1)y$.
When $y\to-\infty$, we have 
$F_{\pm\pm\pm}(y)\sim(\mp\alpha\mp\beta\mp\gamma-1)y$.

We first deal with the case where $\varepsilon A\alpha$, $\varepsilon
B\beta$ and $\varepsilon C\gamma$ are all positive. Then the functions 
$F_{\pm\pm\pm}$ are defined on the whole $\R$. Moreover we have
$\alpha+\beta-\gamma-1<0$ and $-\alpha-\beta+\gamma-1<0$ because
$(\alpha,\beta,\gamma)\in\cK$. Thus by the intermediate value theorem
equation $F_{++-}(y)=0$ has a solution $y\in\R$.

We now deal with the case where at least one of the numbers
$\varepsilon A\alpha$, $\varepsilon B\beta$ and $\varepsilon C\gamma$
is negative (for example $\varepsilon B\beta$). We can assume that
$\frac{\varepsilon B}{2\pi\beta}\leqslant\frac{\varepsilon A}{2\pi\alpha}$ and
$\frac{\varepsilon B}{2\pi\beta}\leqslant\frac{\varepsilon C}{2\pi\gamma}$.
Then the functions $F_{\pm,\pm,\pm}$ are defined for $|y|\geqslant y_0$ where
$y_0=\sqrt{\frac{|B|}{2\pi\beta}}$. We have
$F_{+++}(y_0)=F_{+-+}(y_0)$, $\alpha+\beta+\gamma-1>0$ and
$\alpha-\beta+\gamma-1<0$ (because
$(\alpha,\beta,\gamma)\in\cK$). Thus by the intermediate value theorem
there exists $y\in\R$ such that $F_{+++}(y)=0$ or $F_{+-+}(y)=0$
(depending on the sign of $F_{+++}(y_0)$).
\end{proof}

\begin{rem}
If $\alpha,\beta,\gamma\in(0,1)$, then we have
$\varepsilon=\varepsilon_0$.
\end{rem}

\subsection{The differential equation satisfied by $K_1$ and $K_2$}

This section contains technical results that will be used to construct
trinoids in hyperbolic space (section \ref{trinoids}).
We assume here that the numbers $p$, $q$ and $r$ defined in
lemma \ref{equivalence} satisfy system (\ref{systempqr}), and that
they are real.

\begin{lemma} \label{lemmaeqdiffK}
Set $\hat\Theta\rmd z^2=\rmS_zh-\rmS_z\zeta$ where $h=\frac{K_2}{K_1}$
and $\zeta(z)=\int_0^z\varphi(\tau)\rmd\tau$. Then $K_1$ and $K_2$ are
solutions on $\Sigma$ of the following differential equation:
\begin{equation} \label{eqK}
K''-\frac{\varphi'}{\varphi}K'+\frac{\hat\Theta}{2}K=0.
\end{equation}
\end{lemma}

\begin{proof}
We know that $K_1$ and $K_2$ satisfy (\ref{eqK1K2}) with
$F=\varepsilon\varphi$ (because $(p,q,r)$ is a solution of
(\ref{systempqr})). Then the proof is the same as that of equation 
(\ref{eqk}) in lemma \ref{lemmak1k2}.
\end{proof}

\begin{lemma} \label{lemmahattheta}
The function $\hat\Theta$ extends to a meromorphic function on
$\bar\C$. Moreover we have
$$\hat\Theta=\frac{\Phi}{z^2(z-1)^2}+\frac{\hat\Lambda}{z(z-1)\varphi}
+\frac{2\varphi''}{\varphi}$$ with
$\Phi$ as in (\ref{Phi}) and wheer $\hat\Lambda$ is an affine function.
\end{lemma}

\begin{proof}
By proposition \ref{threelines} there exist $\lambda,\mu\in\C^*$ such
that the functions $k_1=\lambda^{-1}K_1$ and $k_2=\mu K_2$ define a
minimal immersion bounded by some $(D_1,D_2,D_3)\in\cD$, with possibly
a singular point if $\varphi$ has a double
root). Then the Gauss map of this immersion is
$g=\lambda\mu h$, so $\rmS_zh=\rmS_zg$. Thus $\hat\Theta$
extends to a meromorphic function on $\bar\C$ by Schwarz reflection,
and the expression of $\hat\Theta$ follows from lemma \ref{Theta}
(this remains true in the case where $\varphi$ has a double root,
since the order of the pole of $\rmS_zh$ at the root of $\varphi$ is
at most $2$).
\end{proof}

%Since $h(z)\sim\frac{a+\alpha c}{a}z^\alpha$ when $z\to 0$, we have
%$\hat\Theta\sim\frac{1-\alpha^2}{2}z^{-2}$ when $z\to 0$, and the same
%at $1$ and $\infty$ (detailler: dire qu'on peut remplacer h par sa
%composee avec une homographie, cf la grosse preuve).
%
%Montrer que $\hat\Theta$ n'a pas de singularites a part 0,1,infini et
%les zeros de F.

\begin{lemma} \label{lemmahatlambda}
We have 
$$\hat\Lambda(0)=
\varepsilon(p+q+r)((\gamma^2-\beta^2)(2p+q+r)+(1-\alpha^2)(q-r)),$$
$$\hat\Lambda(1)=
\varepsilon(p+q+r)((\alpha^2-\beta^2)(p+q+2r)+(1-\gamma^2)(q-p)),$$
$$\hat\Lambda(1)-\hat\Lambda(0)=
\varepsilon(p+q+r)((\alpha^2-\gamma^2)(p+2q+r)+(1-\beta^2)(r-p)).$$
\end{lemma}

\begin{proof}
We recall that in the neighbourhood of $0$ we have
$$K_j=z^{\frac{1-\alpha}{2}}(1-z)^{\frac{1-\gamma}{2}}
\left((a+bz)w_j^{(0)}+cz(1-z)(w_j^{(0)})'\right)$$ for $j=1,2$, with
$w_j^{(0)}$ as in section \ref{hypergeometricequation}. We have
$$h(z)=z^\alpha\frac{a+\alpha c+\left(b-\alpha c
+\frac{s_{+--}s_{++-}}{1+\alpha}(a+(1+\alpha)c)\right)z+\rmO(z^2)}
{a+\left(b+\frac{s_{---}s_{-+-}}{1-\alpha}(a+c)\right)z+\rmO(z^2)}.$$
The coefficient of the order $-1$ term in $\hat\Theta$ at $0$ is $\hat
s_{-1}=\frac{1-\alpha^2}{\alpha}\frac{h_1}{h_0}$ where
$h(z)=z^\alpha(h_0+h_1z+\rmO(z^2))$. We compute that
\begin{eqnarray*}
\hat s_{-1} & = & \frac{1}{2(p^2-\alpha^2(p+q+r)^2)}\times \\
& & (\alpha^2(\alpha^2-\beta^2+\gamma^2-1)(p+q+r)^2
+(-\alpha^2+5\beta^2-5\gamma^2+1)p^2 \\
& & +(2\alpha^2+2\beta^2-2\gamma^2-2)q^2
+(-2\alpha^2+2\beta^2-2\gamma^2+2)r^2 \\
& & +2(\alpha^2+3\beta^2-3\gamma^2-1)pq
+2(-\alpha^2+3\beta^2-3\gamma^2+1)pr \\
& & +4(\beta^2-\gamma^2)qr).
\end{eqnarray*}
Using that $p^2-\alpha^2(p+q+r)^2=\varepsilon\frac{A\alpha}{2\pi}$ we
get that
\begin{eqnarray*}
\hat s_{-1} & = & \varepsilon\frac{\pi}{A\alpha}\times \\
& & (-(\alpha^2-\beta^2+\gamma^2-1)\varepsilon\frac{A\alpha}{2\pi}
+4(\beta^2-\gamma^2)p^2 \\
& & +(2\alpha^2+2\beta^2-2\gamma^2-2)q^2
+(-2\alpha^2+2\beta^2-2\gamma^2+2)r^2 \\
& & +2(\alpha^2+3\beta^2-3\gamma^2-1)pq
+2(-\alpha^2+3\beta^2-3\gamma^2+1)pr \\
& & +4(\beta^2-\gamma^2)qr).
\end{eqnarray*}
On the other hand we have 
$\hat s_{-1}=-\frac{\hat\Lambda(0)}{\varphi(0)}+\Phi'(0)+2\Phi(0)$, so
using the fact that
$\Phi'(0)+2\Phi(0)=\frac{1-\alpha^2+\beta^2-\gamma^2}{2}$ and
$\varphi(0)=\frac{A\alpha}{2\pi}$
we conclude that 
$$\hat\Lambda(0)=\varepsilon(p+q+r)
((\gamma^2-\beta^2)(2p+q+r)+(1-\alpha^2)(q-r)).$$

In the neighbourhood of $1$ we have 
$$K_j=z^{\frac{1-\alpha}{2}}(1-z)^{\frac{1-\gamma}{2}}
\left((a+bz)w_j^{(1)}+cz(1-z)(w_j^{(1)})'\right)$$ for $j=1,2$, with
$w_j^{(1)}$ as in section \ref{hypergeometricequation}. In the same
way we compute that $$\hat\Lambda(1)=
\varepsilon(p+q+r)((\alpha^2-\beta^2)(p+q+2r)+(1-\gamma^2)(q-p)),$$
and we deduce the expression of
$\hat\Lambda(1)-\hat\Lambda(0)$.
\end{proof}

\section{Application to trinoids in hyperbolic space} \label{trinoids}

\subsection{The cousin relation}

We recall a few facts about the cousin relation between minimal
surfaces in $\R^3$ and Bryant surfaces, {\it i.e.}
constant mean curvature one (CMC-$1$) surfaces 
in hyperbolic space $\h^3$. The
asymptotic boundary of $\h^3$ will be denoted by $\partial_\infty\h^3$.

Let $\cS$ be a simply connected Riemann surface. If $x:\cS\to\R^3$
is a conformal minimal immersion, $\rmI$ and $\rmII$ its first and
second fundamental forms, then
there exists a conformal CMC-$1$ immersion $\tilde{x}:\cS\to\h^3$
whose first and second fundamental forms are
$$\tilde\rmI=\rmI,\quad\widetilde{\rmII}=\rmII+\rmI,$$ and
conversely. The immersions $x$ and $\tilde{x}$ are said to be cousin
immersions. They are unique up to isometries of $\R^3$ and
$\h^3$ respectively.

Bryant proved in \cite{bryant} that if $\tilde{x}$ is such
an immersion (with $\cS$ non necessarily simply connected), then there
exists a holomorphic immersion 
$F:\tilde\cS\to\mathrm{SL}_2(\C)$ where $\tilde\cS$ is the universal
cover of $\cS$ such that $\tilde{x}=FF^*$ and
$\det(F^{-1}\rmd F)=0$, where the model of hyperbolic space is
$$\h^3=\{M\in{\mathcal M}_2(\C)|M^*=M,\tr{M}>0,\det{M}=1\}.$$
Moreover we have
$$F^{-1}\rmd F=\left(\begin{array}{cc}
g & -g^2 \\
1 & -g
\end{array}\right)\omega$$ where $(g,\omega)$ are the Weierstrass data
of the cousin immersion $x$ (see also \cite{umehara}). The map $F$ is
called the Bryant representation of $x$.

The geodesic lines of curvature of $x$ correspond to the geodesic
lines of curvature of $\tilde{x}$, and they lie in planes that are
orthogonal to the surface. The Schwarz reflexion principle
for geodesic lines of curvature also holds for Bryant surfaces. Thus a
planar symmetry of $x$ correspond to a planar symmetry of
$\tilde{x}$. These facts are explained in details in \cite{karcher}
and \cite{toubiana}.

The cousin immersion of the conjugate immersion of $x$ will be called
the conjugate cousin immersion of $x$ and it will be denoted by
$x^\circ$. Thus the Weierstrass data of $x^\circ$ are
$(g^\circ,\omega^\circ)=(g,i\omega)$. Moreover, straight lines of $x$
correspond to geodesic lines of curvature of $x^\circ$ (hence lying in
hyperbolic planes), and symmetries of $x$ with respect to a straight
line correspond to symmetries of $x^\circ$ with respect to a
hyperbolic plane. 

\subsection{Trinoids}

S\'a Earp and Toubiana proved in \cite{toubiana} that an embedded end
of finite total curvature is either asymptotic to the end of a
rotational catenoid
cousin (in which case it is called a catenoidal end) or to a horosphere
(in which case it is called a horospherical end). The catenoidal ends
are the embedded type I ends in the sense of \cite{umehara}; they are
asymptotically rotational surfaces (see \cite{flux}).

\begin{defn}
A Bryant surface is called a trinoid if it has genus zero and three
catenoidal ends.
\end{defn}

Consequently, a trinoid is given by a conformal CMC-1 immersion
defined on $\C\setminus\{0,1\}$. The three ends correspond to $0$, $1$
and $\infty$.

The definitions of a catenoidal end and of a trinoid in \cite{bobenko} are
slightly different from ours. However it turns out that the
definitions of trinoids are equivalent. 

In \cite{collin}, Collin, Hauswirth and Rosenberg proved many results
about properly embedded Bryant surfaces: a properly embedded Bryant surface
of genus $0$ with $1$ end (respectively $2$ ends, $3$ ends) is a
horosphere (respectively an embedded catenoid cousin, an embedded trinoid).

The aim of this section is to construct trinoids by the method of the
conjugate cousin immersion. We first prove the following proposition,
which is a reformulation of lemma 2.4 in \cite{toubiana} and which
will be useful in the sequel.

\begin{prop} \label{holomorphicschwarzian}
Let $\cO$ be a neighbourhood of $0$ in $\C$, $\cO^*=\cO\setminus\{0\}$
and $\widetilde{\cO^*}$ be the universal cover of $\cO^*$. Let
$\mu\in\R^*$ and let $x:\widetilde{\cO^*}\to\h^3$ be a conformal
CMC-$1$ immersion whose Weierstrass data $(g,\omega)$ satisfy
$$g(z)\sim g_0z^\mu,\quad\omega\sim\omega_0z^{-1-\mu}\rmd z$$ when
$z\to 0$ with $g_0,\omega_0\in\C^*$. Let $Q$ be its Hopf differential.

Then $x$ is an embedding of a punctured neighbourhood of $0$ in
$\cO^*$ if and only if the $2$-form $\rmS_zg-2Q$ is holomorphic at
$0$.
\end{prop}

\begin{proof}
We set $$g(z)=z^{\mu}(g_0+g_1z+\rmO(z^2)),\quad
\omega=z^{-1-\mu}(\omega_0+\omega_1z+\rmO(z^2))\rmd z.$$ Then we have
$$Q=z^{-2}(q_{-2}+q_{-1}z+\rmO(1))\rmd z^2,\quad
\rmS_zg=z^{-2}(s_{-2}+s_{-1}z+\rmO(1))\rmd z^2$$
with $q_{-2}=\mu\omega_0g_0$,
$q_{-1}=\mu\omega_1g_0+(1+\mu)\omega_0g_1$,
$s_{-2}=\frac{1-\mu^2}{2}$,
$s_{-1}=\frac{1-\mu^2}{\mu}\frac{g_1}{g_0}$.

We first assume that $\mu>0$. Then,
according to lemma 2.4 in \cite{toubiana}, $x$ is an embedding
of a punctured neighbourhood of $0$ in $\cO^*$ if and only if
\begin{equation} \label{conditiontoubiana}
g_0\omega_0=\frac{1-\mu^2}{4\mu},\quad
(1+\mu)\frac{\omega_1}{\omega_0}
=2\mu\omega_1g_0+2(1+\mu)\omega_0g_1.
\end{equation}
The first condition in (\ref{conditiontoubiana}) is equivalent to
$s_{-2}=2q_{-2}$ (this means that $\rmS_zg-2Q$ has at most a pole of
order $1$ at $0$). If this condition is satisfied, then
since $\frac{q_{-1}}{q_{-2}}=\frac{\omega_{1}}{\omega_{0}}
+\frac{1+\mu}{\mu}\frac{g_1}{g_0}$, the second condition in
(\ref{conditiontoubiana}) is equivalent to
$(1+\mu)\left(\frac{q_{-1}}{q_{-2}}-\frac{s_{-1}}{1-\mu}\right)=2q_{-1},$
and thus to $s_{-1}=2q_{-1}$, which completes the proof in the case
where $\mu>0$.

We now deal with the case where $\mu<0$. The data
$(g^{-1},-g^2\omega)$ define the same CMC-$1$ immersion as
$(g,\omega)$ (see for example \cite{umehara}), and we have
$g^{-1}(z)\sim g_0^{-1}z^{-\mu}$, 
$-g^2\omega\sim-g_0^2\omega_0z^{-1+\mu}\rmd z$, and $\rmS_zg-2Q$ is
unchanged, so it suffices to apply the previous case with $-\mu$
instead of $\mu$.
\end{proof}

\begin{rem}
Umehara and Yamada proved in \cite{umehara} that $\rmS_zG=\rmS_zg-2Q$
where $G$ is the hyperbolic Gauss map of the immersion $x$.
\end{rem}

Let $(D_1,D_2,D_3)\in\cD$,
$(\alpha_0,\gamma_0,\beta_0,-A,-C,-B,\varepsilon_0)=L(D_1,D_2,D_3)$
(see section \ref{geometricconfiguration}), $\alpha\in\alpha_0+2\Z$,
$\beta\in\beta_0+2\Z$, $\gamma\in\gamma_0+2\Z$.
Let $x:\Sigma\to\R^3$ be a minimal immersion bounded by
$(D_1,D_2,D_3)$ or its dual configuration, and
having helicoidal ends of parameters $(A,\alpha)$, 
$(B,\beta)$ and $(C,\gamma)$ at $0$, $1$ and $\infty$ respectively,
corresponding to a solution $(p,q,r)$ of (\ref{systempqr}) where
$\varepsilon=\pm 1$. Let $(g,\omega)$
be its Weierstrass data and $Q=\omega\rmd g$ its Hopf differential.

Then the conjugate cousin immersion $x^\circ:\Sigma\to\h^3$ maps
$(-\infty,0)$, $(0,1)$ and $(1,\infty)$ onto three geodesic lines of
curvature $\cL_1$, $\cL_2$ and $\cL_3$ belonging to three hyperbolic
planes $\cP_1$, $\cP_2$ and $\cP_3$.

\begin{prop} \label{asymptoticboundary}
If $\frac{\alpha^2}{4}>\frac{A\alpha}{2\pi}$,
$\frac{\beta^2}{4}>\frac{B\beta}{2\pi}$ and
$\frac{\gamma^2}{4}>\frac{C\gamma}{2\pi}$, then the asymptotic
boundary of each end of $x^\circ$ consists of one point. 
\end{prop}

\begin{proof}
It suffices to prove that the asymptotic boundary of the end at $0$
consists of one point. The Hopf differential
$Q^\circ=\omega^\circ\rmd g^\circ$ satisfies
$Q^\circ\sim q_{-2}z^{-2}\rmd z$ at $0$ with
$q_{-2}=-\frac{A\alpha}{2\pi}$. We proceeding as in the proof of lemma
2.4 in \cite{toubiana}: since the indicial equations
$\tau^2+\alpha\tau-q_{-2}$ and $\upsilon^2-\alpha\upsilon-q_{-2}$ have
a positive discriminant $\Delta=\alpha^2+4q_{-2}$ (because of the
hypothesis), we prove that, up to an isometry of $\h^3$, the 
Bryant representation of $x$ is
$F=\left( \begin{array}{cc}
z^\tau A_1(z) & z^\upsilon B_1(z) \\
z^\tau C_1(z) & z^\upsilon D_1(z) 
\end{array} \right)$ where $\tau=\frac{-\sqrt\Delta-\alpha}{2}$,
$\upsilon=\frac{-\sqrt\Delta+\alpha}{2}$ and where $A_1$, $B_1$,
$C_1$ and $D_1$ 
are holomorphic functions in a neighbourhood of $0$ in 
$\{\im z\geqslant 0\}$ that do not vanish at $0$. Using the
identification of $\h^3$ with the upper half-space model
$\{(y_1,y_2,y_3)\in\R^3|y_3>0\}$ described in \cite{toubiana}, we get
$$y_1+iy_2=\frac{|z|^{2\tau}\overline{A_1}C_1
+|z|^{2\upsilon}\overline{B_1}D_1}
{|z|^{2\tau}|A_1|^2+|z|^{2\upsilon}|B_1|^2},\quad
y_3=\frac{1}{|z|^{2\tau}|A_1|^2+|z|^{2\upsilon}|B_1|^2}$$
(this is formula (1.2) of \cite{toubiana}). Thus we have $y_3\to 0$
when $z\to 0$ (since $\tau$ or $\upsilon$ is negative), and
$y_1+iy_2\to\frac{C_1(0)}{A_1(0)}$ or
$y_1+iy_2\to\frac{D_1(0)}{B_1(0)}$ (depending on the sign of
$\alpha$). This proves the assertion.
\end{proof}

From now on we assume that $\frac{\alpha^2}{4}>\frac{A\alpha}{2\pi}$,
$\frac{\beta^2}{4}>\frac{B\beta}{2\pi}$ and
$\frac{\gamma^2}{4}>\frac{C\gamma}{2\pi}$.
Thus the lines $\cL_1$, $\cL_2$ and $\cL_3$ are pairwise
concurrent in $\partial_\infty\h^3$ at the asymptotic boundary
points. Applying Schwarz reflections with 
respect to these planes and repeating the process with respect to the
new planes infinitely many times, we get a conformal CMC-1 immersion
$x^\circ:\widetilde{\C\setminus\{0,1\}}\to\h^3$ where
$\widetilde{\C\setminus\{0,1\}}$ is the universal cover of
$\C\setminus\{0,1\}$. This immersion $x^\circ$ is well-defined on
$\C\setminus\{0,1\}$ if and only if the planes $\cP_1$, $\cP_2$ and
$\cP_3$ are equal.

\begin{prop} \label{necessarycondition}
If immersion $x^\circ$ gives a trinoid by Schwarz reflection, then we
have 
\begin{equation} \label{conditiontrinoid}
\frac{A\alpha}{2\pi}=\frac{\alpha^2-1}{4},\quad
\frac{B\beta}{2\pi}=\frac{\beta^2-1}{4},\quad
\frac{C\gamma}{2\pi}=\frac{\gamma^2-1}{4}.
\end{equation}
\end{prop}

\begin{proof}
Assume that $x^\circ$ gives a trinoid by Schwarz reflection. Then
$x^\circ$ is well-defined on $\C\setminus\{0,1\}$, and its ends are
embedded. Let $Q^\circ$ be its Hopf differential. We have $Q^\circ=iQ$.

Its Weierstrass data satisfy $g^\circ(z)=g(z)\sim g_0z^\alpha$ and
$\omega^\circ=i\omega\sim\omega_0z^{-1-\alpha}\rmd z$ when $z\to 0$,
with $g_0,\omega_0\in\C^*$. Then by proposition
\ref{holomorphicschwarzian} the form $\rmS_zg^\circ-2Q^\circ$ is
holomorphic at $0$. In particular the order $-2$ term vanishes, that
is $\frac{1-\alpha^2}{2}+2\frac{A\alpha}{2\pi}=0$.

The other two identities are obtained in the same way for the ends at
$\infty$ and $1$.
\end{proof}

A computation gives the following result.

\begin{lemma} \label{lemmasystempqrtrinoid}
The complex solutions of system 
\begin{eqnarray} \label{systempqrtrinoid}
\left\{\begin{array}{ccc}
p^2-\alpha^2(p+q+r)^2 & = & \frac{\alpha^2-1}{4} \\
q^2-\beta^2(p+q+r)^2 & = & \frac{\beta^2-1}{4} \\
r^2-\gamma^2(p+q+r)^2 & = & \frac{\gamma^2-1}{4}
\end{array}\right.
\end{eqnarray}
are $(-\frac{i}{2},\frac{i}{2},\frac{i}{2})$,
$(\frac{i}{2},-\frac{i}{2},\frac{i}{2})$,
$(\frac{i}{2},\frac{i}{2},-\frac{i}{2})$,
$(U\delta,V\delta,W\delta)$ and their opposites, where
$$U=-3\alpha^4+2(1+\beta^2+\gamma^2)\alpha^2+\beta^4+\gamma^4
-2(\beta^2+\gamma^2+\beta^2\gamma^2)+1,$$
$$V=-3\beta^4+2(1+\alpha^2+\gamma^2)\beta^2+\alpha^4+\gamma^4
-2(\alpha^2+\gamma^2+\alpha^2\gamma^2)+1,$$
$$W=-3\gamma^4+2(1+\alpha^2+\beta^2)\gamma^2+\alpha^4+\beta^4
-2(\alpha^2+\beta^2+\alpha^2\beta^2)+1,$$ and where $\delta$ is a
complex square root of $-\frac{1}{4\Pi}$ with $\Pi$ defined by (\ref{Pi}).

The complex solutions of system 
\begin{eqnarray} \label{systempqrtrinoidnegative}
\left\{\begin{array}{ccc}
p^2-\alpha^2(p+q+r)^2 & = & \frac{1-\alpha^2}{4} \\
q^2-\beta^2(p+q+r)^2 & = & \frac{1-\beta^2}{4} \\
r^2-\gamma^2(p+q+r)^2 & = & \frac{1-\gamma^2}{4}
\end{array}\right.
\end{eqnarray}
are $(\frac{1}{2},\frac{1}{2},-\frac{1}{2})$, 
$(\frac{1}{2},-\frac{1}{2},\frac{1}{2})$,
$(-\frac{1}{2},\frac{1}{2},\frac{1}{2})$, 
$(iU\delta,iV\delta,iW\delta)$ and their opposites.
\end{lemma}

\begin{rem} \label{distinctsolutions}
These solutions are distinct if and only if
$1-\alpha^2-\beta^2+\gamma^2\neq 0$, $1-\alpha^2+\beta^2-\gamma^2\neq
0$ and $1+\alpha^2-\beta^2-\gamma^2\neq 0$.
\end{rem}

\begin{rem} \label{discrPhi}
We compute that $U+V+W$ is $-4$ times the discriminant of $\Phi$.
\end{rem}

\begin{prop} \label{goodsolution}
The immersion $x^\circ$ gives a trinoid by Schwarz reflection if and
only if (\ref{conditiontrinoid}) holds, $\varphi$ has no double root, and
$(p,q,r)=\pm(U\delta,V\delta,W\delta)$ in the case
where $\delta\in\R$ or $(p,q,r)=\pm(iU\delta,iV\delta,iW\delta)$ in the case
where $\delta\in i\R$, where $U$, $V$, $W$ and $\delta$ are as in
lemma \ref{lemmasystempqrtrinoid}.
\end{prop}

\begin{proof}
We have $$\rmS_zg^\circ-2Q^\circ=
\left(\frac{\Phi+2\varphi}{z^2(z-1)^2}+\frac{\Lambda}{\varphi z(z-1)}
+\frac{2\varphi''}{\varphi}\right)
\rmd z^2+\rmS_z\zeta$$
with $\zeta(z)=\int_0^z\varphi(\tau)\rmd\tau$, $\Phi$ as in
(\ref{Phi}) and $\Lambda$ as $\hat\Lambda$ in lemma
\ref{lemmahatlambda}. 
By proposition \ref{holomorphicschwarzian}, $x^\circ$ gives a trinoid
if and only if $\rmS_zg^\circ-Q^\circ$ is holomorphic at $0$, $1$ and
$\infty$. This holds if and only if $\Phi=-2\varphi$ ({\it i.e.} if
(\ref{conditiontrinoid}) holds, {\it i.e.} if $(p,q,r)$ is a real
solution of (\ref{systempqrtrinoid}) or
(\ref{systempqrtrinoidnegative})) and $\Lambda=0$ (we recall that
$\frac{2\varphi''}{\varphi}\rmd z^2+\rmS_z\zeta$ is holomorphic at $\infty$).

By lemma \ref{lemmahatlambda} we have $\Lambda=0$ if and only if
$p+q+r=0$ or
\begin{eqnarray} \label{systemlambda}
\left\{ \begin{array}{ccc}
(\gamma^2-\beta^2)(2p+q+r)+(1-\alpha^2)(q-r) & = & 0 \\
(\alpha^2-\beta^2)(p+q+2r)+(1-\gamma^2)(q-p) & = & 0.
\end{array} \right.
\end{eqnarray}

We notice that $(U,V,W)\neq(0,0,0)$ (since $(U\delta,V\delta,W\delta)$
is solution of (\ref{systempqrtrinoid})). Thus the set of the
solutions of (\ref{systemlambda}) is the complex line $\C(U,V,W)$. Hence
the only solutions of (\ref{systemlambda}) that are also solutions of
(\ref{systempqrtrinoid}) or (\ref{systempqrtrinoidnegative}) are
$(U\delta,V\delta,W\delta)$, $(iU\delta,iV\delta,iW\delta)$ and their
opposites. 

There is a solution $(p,q,r)$ of (\ref{systempqrtrinoid}) or
(\ref{systempqrtrinoidnegative}) satisfying $p+q+r=0$ if and only if
$U+V+W=0$, {\it i.e.} if and only if $\varphi$ has a double real root
$a_1$ (by remark \ref{discrPhi} and since $\Phi=-2\varphi$). In this
case these solutions are again $(U\delta,V\delta,W\delta)$,
$(iU\delta,iV\delta,iW\delta)$ and their opposites.

Hence, in the case where $\varphi$ has a double root $a_1$, the
solutions $(p,q,r)$ satisfy $p+q+r=0$, which implies $\Lambda=0$ by lemma
\ref{lemmahatlambda}. Thus the exponents of equation (\ref{eqk}) at
$a_1$ are $1$ and $2$. This implies that the spinors $k_1$ and $k_2$
associated to $x$ both vanish at $a_1$, and so $x$ and $x^\circ$ do
have a singular point at $a_1$. 

Moreover, the numbers $p$, $q$ and $r$ are required to be real, so the
proof is complete.
\end{proof}

\begin{rem}
In the case where $\varphi$ has a double root, the result still holds
except that the immersion giving the trinoid has a singularity at the
root of $\varphi$.
\end{rem}

\begin{thm} \label{theoremtrinoids}
Let $\mu_0$, $\mu_1$ and $\mu_\infty$ be three positive non-integer real
numbers. Assume that
\begin{equation} \label{conditiongrowths1}
(|[\mu_0]|,|[\mu_1]|,|[\mu_\infty]|)\in\cK,
\end{equation}
where $[r]$ denotes the unique number in $(-1,1]$ such that
$r-[r]\in2\Z$ (the set $\cK$ is defined in proposition
\ref{bijection}), and that 
\begin{equation} \label{conditiongrowths2}
\mu_0^4+\mu_1^4+\mu_\infty^4
-2\mu_0^2\mu_1^2-2\mu_0^2\mu_\infty^2-2\mu_1^2\mu_\infty^2
+2\mu_0^2+2\mu_1^2+2\mu_\infty^2-3\neq 0.
\end{equation}
Then there exists a trinoid $\cT_{\mu_0,\mu_1,\mu_\infty}$ whose
ends are of growths $1-\mu_0$, $1-\mu_1$ and $1-\mu_\infty$ and having
a symmetry plane. 

The ends of $\cT_{\mu_0,\mu_1,\mu_\infty}$ have distinct asymptotic
boundary points if and only
if $1-\mu_0^2-\mu_1^2+\mu_\infty^2\neq 0$,
$1-\mu_0^2+\mu_1^2-\mu_\infty^2\neq 0$ and
$1+\mu_0^2-\mu_1^2-\mu_\infty^2\neq 0$. More precisely, the ends of
growths $1-\mu_0$ and $1-\mu_1$ (respectively $1-\mu_0$ and
$1-\mu_\infty$, $1-\mu_1$ and $1-\mu_\infty$) have distinct asymptotic
boundary points if and only if $1-\mu_0^2-\mu_1^2+\mu_\infty^2\neq 0$
(respectively $1-\mu_0^2+\mu_1^2-\mu_\infty^2\neq 0$,
$1+\mu_0^2-\mu_1^2-\mu_\infty^2\neq 0$). In particular the three ends
cannot have the same asymptotic boundary point.
\end{thm} 

\begin{proof}
We set $\alpha=\mu_0$ if $\mu_0\in|[\mu_0]|+2\Z$ and
$\alpha=-\mu_0$ if $-\mu_0\in|[\mu_0]|+2\Z$ (in order to be compatible
with the conventions of section \ref{hopfandspinors}). In the same way we set
$\beta=\pm\mu_\infty$ and $\gamma=\pm\mu_1$. Let
$\varepsilon_0\in\{1,-1\}$. By proposition \ref{bijection} there
exists a triple $(D_1,D_2,D_3)$ such that
$L(D_1,D_2,D_3)=(|[\mu_0]|,|[\mu_1]|,|[\mu_\infty]|,-A,-C,-B,\varepsilon_0)$
with $\frac{A\alpha}{2\pi}=\frac{\alpha^2-1}{4}$,
$\frac{B\beta}{2\pi}=\frac{\beta^2-1}{4}$,
$\frac{C\gamma}{2\pi}=\frac{\gamma^2-1}{4}$.
Then by (\ref{conditiongrowths2}) the corresponding $\varphi$ has no
double root, and so by proposition \ref{goodsolution} there exists a
minimal immersion
$x:\Sigma\to\R^3$ bounded by $(D_1,D_2,D_3)$ or its dual configuration
whose conjugate cousin $x^\circ$ gives a trinoid by Schwarz reflection
; moreover the growths of the ends of this trinoid are $1-\mu_0$,
$1-\mu_1$ and $1-\mu_\infty$ respectively (it suffices to consider the
coefficient of the order $-2$ term of $Q^\circ$ at each end). This
proves the existence of the trinoid $\cT_{\mu_0,\mu_1,\mu_\infty}$ (it
has a symmetry plane by construction).

Up to an isometry of $\h^3$, the hyperbolic Gauss map of $x^\circ$ is
$G^\circ(z)=z+\frac{(a_1-a_2)^2}{2(2z-a_1-a_2)}$ where $a_1$ and $a_2$
are the roots of $\varphi$ (see \cite{period}, example 4.4; we notice
that $\rmS_zG^\circ=\rmS_z\zeta+\frac{2\varphi''}{\varphi}\rmd z^2$).
Moreover, the limit of the hyperbolic Gauss map at a catenoidal end is
the asymptotic boundary point of the end (see \cite{toubiana}). Thus,
to compare the asymptotic boundary points of the ends, it suffices to
compare $G^\circ(0)$, $G^\circ(1)$ and $G^\circ(\infty)$.

We have $G^\circ(0)=-\frac{(a_1-a_2)^2}{2(a_1+a_2)}$,
$G^\circ(1)=1+\frac{(a_1-a_2)^2}{2(2-a_1-a_2)}$ and
$G^\circ(\infty)=\infty$. Thus we have $G^\circ(0)=G^\circ(1)$ if and
only if $a_1+a_2=2a_1a_2$, {\it i.e.}
$1-\mu_0^2-\mu_1^2+\mu_\infty^2=0$;
we have $G^\circ(0)=G^\circ(\infty)$ if
and only if $a_1+a_2=0$, {\it i.e.} $1-\mu_0^2+\mu_1^2-\mu_\infty^2=0$
; we have $G^\circ(1)=G^\circ(\infty)$ if and only if $a_1+a_2=2$,
{\it i.e.} $1+\mu_0^2-\mu_1^2-\mu_\infty^2=0$. This also implies that
we never have $G^\circ(0)=G^\circ(1)=G^\circ(\infty)$.
\end{proof}

\begin{rem}
If (\ref{conditiongrowths1}) holds but (\ref{conditiongrowths2})
does not hold, then there exists a ``trinoid''
$\cT_{\mu_0,\mu_1,\mu_\infty}$ with one singular point.
\end{rem}

\begin{cor}
If $(\mu_0,\mu_1,\mu_\infty)\in\cK$ and 
$\mu_0,\mu_1,\mu_\infty\in(0,1)$ ({\it i.e.} if the
growths are positive), then the trinoid $\cT_{\mu_0,\mu_1,\mu_\infty}$
exists and its ends have distinct asymptotic boundary points.
\end{cor}

\begin{proof}
In this case we have
$\mu_j=|[\mu_j]|$ for $j=0,1,\infty$, and so
$(\mu_0,\mu_1,\mu_\infty)\in\cK$. This implies that
$\mu_\infty>1-\mu_0-\mu_1$ and $\mu_\infty>-1+\mu_0+\mu_1$, and so
$1-\mu_0^2-\mu_1^2+\mu_\infty^2>1-\mu_0^2-\mu_1^2+(1-\mu_0-\mu_1)^2 
=2(1-\mu_0)(1-\mu_1)>0$. In the same way we have 
$1-\mu_0^2+\mu_1^2-\mu_\infty^2\neq 0$ and
$1+\mu_0^2-\mu_1^2-\mu_\infty^2\neq 0$. Thus the asymptotic boundary points
of the ends are distinct. 

Set $d(\mu_0,\mu_1,\mu_\infty)=\mu_0^4+\mu_1^4+\mu_\infty^4
-2\mu_0^2\mu_1^2-2\mu_0^2\mu_\infty^2-2\mu_1^2\mu_\infty^2
+2\mu_0^2+2\mu_1^2+2\mu_\infty^2-3$. The derivative of $d$ with respect
to $\mu_\infty^2$ is equal to $2(\mu_\infty^2-\mu_0^2-\mu_1^2+1)$,
which was proven to be positive. Without loss of generality we can
assume that $\mu_0\geqslant\mu_1$. We have $\mu_\infty<1-\mu_0+\mu_1$,
and so $d(\mu_0,\mu_1,\mu_\infty)<d(\mu_0,\mu_1,1-\mu_0+\mu_1)
=8(\mu_1+1)(\mu_0-1)(\mu_0-\mu_1)\leqslant 0$. Thus (\ref{conditiongrowths2})
is satisfied, and the trinoid has no singularity.
\end{proof}

Irreducible trinoids are classified by theorem 2.6 of
\cite{metrics}. They correspond to trinoids with
non-integer growth ends. Theorem 2.6 of
\cite{metrics} states that there exists a trinoid
$\cT_{\mu_0,\mu_1,\mu_\infty}$ (without assuming that it has a
symmetry plane) if and only if (\ref{conditiongrowths2}) holds and
\begin{equation} \label{conditionumehara}
\cos^2(\pi\mu_0)+\cos^2(\pi\mu_1)+\cos^2(\pi\mu_\infty)
+2\cos(\pi\mu_0)\cos(\pi\mu_1)\cos(\pi\mu_\infty)<1,
\end{equation}
and in this case this trinoid is unique (in this theorem, the
$\beta_j$ ($j=1,2,3$) correspond to our $\mu_j-1$ ($j=0,1,\infty$),
the $c_j$ to our $\frac{1-\mu_j^2}{2}$).
Irreducible trinoids are also classified in \cite{bobenko}; it is also
proved that (\ref{conditionumehara}) is equivalent to
(\ref{conditiongrowths1})
(in \cite{bobenko} the $\Delta_j$ ($j=1,2,3$) correspond to our
$\frac{|[\mu_j]|}{2}$ ($j=0,1,\infty$)). Pictures of trinoids can be
found in \cite{bobenko} (see also \cite{period}, example 4.4).

\begin{prop} \label{proofrosenberg}
If the asymptotic boundary points of the ends of $x^\circ$ are
distinct and (\ref{conditiontrinoid}) holds then $x^\circ$ gives a
trinoid by Schwarz reflection.
\end{prop}

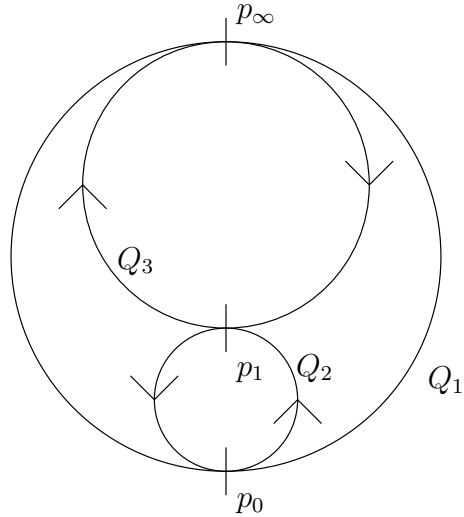
\begin{figure}[htbp]
\begin{center}
\input{circles.pstex_t}
\caption{The asymptotic boundaries of the boundary planes; there is
no orientation of $Q_1$ compatible with those of $Q_2$ and $Q_3$.}
\label{3circles}
\end{center}
\end{figure}

\begin{proof}
For $j=1,2,3$, the oriented curve $\cL_j$ and the mean curvature
vector of $x^\circ$ on $\cL_j$ induce an orientation of the plane $\cP_j$.
Denote by $p_0$, $p_1$ and $p_\infty$ the asymptotic boundaries of the
ends of $x^\circ$ at $0$, $1$ and $\infty$ respectively. 
Then $\cL_1$ goes from $p_\infty$ to $p_0$, $\cL_2$
goes from $p_0$ to $p_1$, and $\cL_3$ goes from $p_1$ to $p_\infty$.
Denote by $Q_1$, $Q_2$ and
$Q_3$ the asymptotic boundaries of $\cP_1$, $\cP_2$ and $\cP_3$. These are
great circles in $\bar\C$. They are given the orientation induced by
$\cP_1$, $\cP_2$ and $\cP_3$ respectively. We have 
$p_0\in Q_1\cap Q_2$, $p_1\in Q_2\cap Q_3$ and 
$p_\infty\in Q_3\cap Q_1$. Moreover, the circles are pairwise tangent at these
points, and their orientations at these points are compatible (since
the boundary lines have turned of an angle $\pi$ at each end, because
of (\ref{conditiontrinoid})). Since $p_0$, $p_1$ and $p_\infty$ are
distinct, this is
not possible unless the three circles are equal (see figure
\ref{3circles}). Consequently, the planes $\cP_1$, $\cP_2$ and $\cP_3$ are
equal, and doing the Schwarz reflection of $x^\circ$ with respect to
this plane gives a trinoid.
\end{proof}

We now describe what the immersion $x^\circ$ looks like in the case
where $(p,q,r)=(\frac{1}{2},\frac{1}{2},-\frac{1}{2})$, 
$(p,q,r)=(\frac{1}{2},-\frac{1}{2},\frac{1}{2})$ or
$(p,q,r)=(-\frac{1}{2},\frac{1}{2},\frac{1}{2})$. We first notice that
$(iU\delta,iV\delta,iW\delta)=(\frac{1}{2},\frac{1}{2},-\frac{1}{2})$,
$(iU\delta,iV\delta,iW\delta)=(\frac{1}{2},-\frac{1}{2},\frac{1}{2})$
and 
$(iU\delta,iV\delta,iW\delta)=(-\frac{1}{2},\frac{1}{2},\frac{1}{2})$
if and only if $1-\alpha^2-\beta^2+\gamma^2=0$,
$1-\alpha^2+\beta^2-\gamma^2=0$ and $1+\alpha^2-\beta^2-\gamma^2=0$
respectively. Henceforth we assume that none of these conditions are
satisfied. 

If $(p,q,r)=(\frac{1}{2},\frac{1}{2},-\frac{1}{2})$, then by lemma
\ref{lemmahatlambda} we have $\Lambda(0)\neq 0$, $\Lambda(1)=0$ and
$\Lambda(1)-\Lambda(0)\neq 0$, so $S_zg^\circ-2Q^\circ$ has poles of
order $1$ at $0$ and $\infty$ and is holomorphic at $1$. This means
that the end at $1$ is embedded, but the ends at $0$ and $\infty$ are
not: the planes $\cP_2$ and $\cP_3$ are identical, and the plane
$\cP_1$ is tangent at infinity to $\cP_2$ but different. Applying
Schwarz reflection infinitely many times, we get a surface that is
invariant by the parabolic isometry generated by the reflections with
respect to $\cP_1$ and $\cP_2$. The ends at $0$ and $\infty$ are not
annular ends since $x^\circ$ is not single-valued at $0$ and $\infty$.

If $(p,q,r)=(\frac{1}{2},-\frac{1}{2},\frac{1}{2})$, then we have 
$\Lambda(0)\neq 0$, $\Lambda(1)\neq 0$ and
$\Lambda(1)-\Lambda(0)=0$, so similarly the end at $\infty$ is
embedded and the ends at $0$ and $1$ are not (and they are not
annular ends).

If $(p,q,r)=(-\frac{1}{2},\frac{1}{2},\frac{1}{2})$, then we have 
$\Lambda(0)=0$, $\Lambda(1)\neq 0$ and
$\Lambda(1)-\Lambda(0)\neq 0$, so the end at $0$ is
embedded and the ends at $1$ and $\infty$ are not (and they are not
annular ends).

Then by propostion \ref{proofrosenberg} the asymptotic boundary points
of these immersions are not pairwise distinct (otherwise these
immersions would give trinoids by Schwarz reflection).

\nocite{*}

\bibliographystyle{alpha}
\bibliography{minimal}

\end{document}

%% file: pospos.pstex_t
\begin{picture}(0,0)%
\includegraphics{pospos.pstex}%
\end{picture}%
\setlength{\unitlength}{2763sp}%
\begingroup\makeatletter\ifx\SetFigFont\undefined%
\gdef\SetFigFont#1#2#3#4#5{%
  \reset@font\fontsize{#1}{#2pt}%
  \fontfamily{#3}\fontseries{#4}\fontshape{#5}%
  \selectfont}%
\fi\endgroup%
\begin{picture}(5424,2424)(4189,-3673)
\put(8476,-2386){\makebox(0,0)[lb]{\smash{\SetFigFont{8}{9.6}{\rmdefault}{\mddefault}{\updefault}{$\pi\alpha$}%
}}}
\put(6826,-2236){\makebox(0,0)[lb]{\smash{\SetFigFont{8}{9.6}{\rmdefault}{\mddefault}{\updefault}{$D_1$}%
}}}
\put(7126,-3061){\makebox(0,0)[lb]{\smash{\SetFigFont{8}{9.6}{\rmdefault}{\mddefault}{\updefault}{$D_2$}%
}}}
\end{picture}

%% file: posneg.pstex_t
\begin{picture}(0,0)%
\includegraphics{posneg.pstex}%
\end{picture}%
\setlength{\unitlength}{2763sp}%
\begingroup\makeatletter\ifx\SetFigFont\undefined%
\gdef\SetFigFont#1#2#3#4#5{%
  \reset@font\fontsize{#1}{#2pt}%
  \fontfamily{#3}\fontseries{#4}\fontshape{#5}%
  \selectfont}%
\fi\endgroup%
\begin{picture}(5424,2424)(4189,-3673)
\put(7126,-3061){\makebox(0,0)[lb]{\smash{\SetFigFont{8}{9.6}{\rmdefault}{\mddefault}{\updefault}{$D_2$}%
}}}
\put(8476,-2311){\makebox(0,0)[lb]{\smash{\SetFigFont{8}{9.6}{\rmdefault}{\mddefault}{\updefault}{$\pi\alpha$}%
}}}
\put(6826,-2236){\makebox(0,0)[lb]{\smash{\SetFigFont{8}{9.6}{\rmdefault}{\mddefault}{\updefault}{$D_1$}%
}}}
\end{picture}

%% file: negpos.pstex_t
\begin{picture}(0,0)%
\includegraphics{negpos.pstex}%
\end{picture}%
\setlength{\unitlength}{2763sp}%
\begingroup\makeatletter\ifx\SetFigFont\undefined%
\gdef\SetFigFont#1#2#3#4#5{%
  \reset@font\fontsize{#1}{#2pt}%
  \fontfamily{#3}\fontseries{#4}\fontshape{#5}%
  \selectfont}%
\fi\endgroup%
\begin{picture}(5424,3575)(4189,-4241)
\put(4576,-1711){\makebox(0,0)[lb]{\smash{\SetFigFont{8}{9.6}{\rmdefault}{\mddefault}{\updefault}{$\pi\alpha$}%
}}}
\put(6826,-2986){\makebox(0,0)[lb]{\smash{\SetFigFont{8}{9.6}{\rmdefault}{\mddefault}{\updefault}{$D_2$}%
}}}
\put(6901,-2236){\makebox(0,0)[lb]{\smash{\SetFigFont{8}{9.6}{\rmdefault}{\mddefault}{\updefault}{$D_1$}%
}}}
\end{picture}

%% file: negneg.pstex_t
\begin{picture}(0,0)%
\includegraphics{negneg.pstex}%
\end{picture}%
\setlength{\unitlength}{2763sp}%
\begingroup\makeatletter\ifx\SetFigFont\undefined%
\gdef\SetFigFont#1#2#3#4#5{%
  \reset@font\fontsize{#1}{#2pt}%
  \fontfamily{#3}\fontseries{#4}\fontshape{#5}%
  \selectfont}%
\fi\endgroup%
\begin{picture}(5424,3575)(4189,-4241)
\put(4651,-1711){\makebox(0,0)[lb]{\smash{\SetFigFont{8}{9.6}{\rmdefault}{\mddefault}{\updefault}{$\pi\alpha$}%
}}}
\put(6901,-2236){\makebox(0,0)[lb]{\smash{\SetFigFont{8}{9.6}{\rmdefault}{\mddefault}{\updefault}{$D_1$}%
}}}
\put(6826,-2986){\makebox(0,0)[lb]{\smash{\SetFigFont{8}{9.6}{\rmdefault}{\mddefault}{\updefault}{$D_2$}%
}}}
\end{picture}

%% file: 3lines.pstex_t
\begin{picture}(0,0)%
\includegraphics{3lines.pstex}%
\end{picture}%
\setlength{\unitlength}{2763sp}%
\begingroup\makeatletter\ifx\SetFigFont\undefined%
\gdef\SetFigFont#1#2#3#4#5{%
  \reset@font\fontsize{#1}{#2pt}%
  \fontfamily{#3}\fontseries{#4}\fontshape{#5}%
  \selectfont}%
\fi\endgroup%
\begin{picture}(6024,2995)(2389,-4544)
\put(7876,-2536){\makebox(0,0)[lb]{\smash{\SetFigFont{8}{9.6}{\rmdefault}{\mddefault}{\updefault}{$\pi\alpha_0$}%
}}}
\put(4351,-4486){\makebox(0,0)[lb]{\smash{\SetFigFont{8}{9.6}{\rmdefault}{\mddefault}{\updefault}{$\pi\beta_0$}%
}}}
\put(5401,-3436){\makebox(0,0)[lb]{\smash{\SetFigFont{8}{9.6}{\rmdefault}{\mddefault}{\updefault}{$D_1$}%
}}}
\put(5101,-2611){\makebox(0,0)[lb]{\smash{\SetFigFont{8}{9.6}{\rmdefault}{\mddefault}{\updefault}{$D_2$}%
}}}
\put(3976,-3286){\makebox(0,0)[lb]{\smash{\SetFigFont{8}{9.6}{\rmdefault}{\mddefault}{\updefault}{$D_3$}%
}}}
\put(2626,-2461){\makebox(0,0)[lb]{\smash{\SetFigFont{8}{9.6}{\rmdefault}{\mddefault}{\updefault}{$\pi\gamma_0$}%
}}}
\end{picture}

%% file: dual.pstex_t
\begin{picture}(0,0)%
\includegraphics{dual.pstex}%
\end{picture}%
\setlength{\unitlength}{2763sp}%
\begingroup\makeatletter\ifx\SetFigFont\undefined%
\gdef\SetFigFont#1#2#3#4#5{%
  \reset@font\fontsize{#1}{#2pt}%
  \fontfamily{#3}\fontseries{#4}\fontshape{#5}%
  \selectfont}%
\fi\endgroup%
\begin{picture}(6099,3670)(4189,-4319)
\put(9601,-2011){\makebox(0,0)[lb]{\smash{\SetFigFont{8}{9.6}{\rmdefault}{\mddefault}{\updefault}{$\pi\alpha_0$}%
}}}
\put(6901,-4261){\makebox(0,0)[lb]{\smash{\SetFigFont{8}{9.6}{\rmdefault}{\mddefault}{\updefault}{$\pi\beta_0$}%
}}}
\put(5101,-1636){\makebox(0,0)[lb]{\smash{\SetFigFont{8}{9.6}{\rmdefault}{\mddefault}{\updefault}{$\pi\gamma_0$}%
}}}
\put(7201,-3061){\makebox(0,0)[lb]{\smash{\SetFigFont{8}{9.6}{\rmdefault}{\mddefault}{\updefault}{$D_2$}%
}}}
\put(6676,-2311){\makebox(0,0)[lb]{\smash{\SetFigFont{8}{9.6}{\rmdefault}{\mddefault}{\updefault}{$D_1$}%
}}}
\put(8476,-2386){\makebox(0,0)[lb]{\smash{\SetFigFont{8}{9.6}{\rmdefault}{\mddefault}{\updefault}{$D_3$}%
}}}
\end{picture}

%% file: circles.pstex_t
\begin{picture}(0,0)%
\includegraphics{circles.pstex}%
\end{picture}%
\setlength{\unitlength}{3947sp}%
\begingroup\makeatletter\ifx\SetFigFont\undefined%
\gdef\SetFigFont#1#2#3#4#5{%
  \reset@font\fontsize{#1}{#2pt}%
  \fontfamily{#3}\fontseries{#4}\fontshape{#5}%
  \selectfont}%
\fi\endgroup%
\begin{picture}(2716,3289)(3143,-3044)
\put(4576, 89){\makebox(0,0)[lb]{\smash{\SetFigFont{12}{14.4}{\rmdefault}{\mddefault}{\updefault}{$p_\infty$}%
}}}
\put(4576,-2986){\makebox(0,0)[lb]{\smash{\SetFigFont{12}{14.4}{\rmdefault}{\mddefault}{\updefault}{$p_0$}%
}}}
\put(5776,-2236){\makebox(0,0)[lb]{\smash{\SetFigFont{12}{14.4}{\rmdefault}{\mddefault}{\updefault}{$Q_1$}%
}}}
\put(4951,-2161){\makebox(0,0)[lb]{\smash{\SetFigFont{12}{14.4}{\rmdefault}{\mddefault}{\updefault}{$Q_2$}%
}}}
\put(3826,-1486){\makebox(0,0)[lb]{\smash{\SetFigFont{12}{14.4}{\rmdefault}{\mddefault}{\updefault}{$Q_3$}%
}}}
\put(4576,-2161){\makebox(0,0)[lb]{\smash{\SetFigFont{12}{14.4}{\rmdefault}{\mddefault}{\updefault}{$p_1$}%
}}}
\end{picture}